\documentclass[12pt, regno]{amsart}
\usepackage[main=english]{babel}
\usepackage[utf8]{inputenc}
\usepackage{amssymb}
\usepackage{amsmath}
\usepackage{amsthm}
\usepackage{setspace}
\usepackage{thmtools}
\usepackage{mathtools}
\usepackage{bbold}
\usepackage{thm-restate}
\usepackage{etoolbox}
\usepackage{fmtcount}
\usepackage{mathrsfs}
\usepackage{fullpage} %fit more on a page without having to set margins manually
\usepackage{colonequals} %Correct :=
\usepackage{enumerate, enumitem}
\usepackage[dvipsnames]{xcolor}
\usepackage{hyperref}
\usepackage{yhmath}
\hypersetup{hypertexnames=false, colorlinks=true, linkcolor=cyan, citecolor=magenta, pdfpagemode=UseNone, pdfstartview={XYZ null null 1.00}}
\usepackage{tikz-cd}
\usetikzlibrary{calc}
\usepackage[capitalise]{cleveref}
\usepackage{indentfirst}

\newlength{\bibitemsep}
\newlength{\bibparskip}
\let\oldthebibliography\thebibliography
\renewcommand\thebibliography[1]{%
  \oldthebibliography{#1}%
  \setlength{\parskip}{\bibitemsep}%
  \setlength{\itemsep}{\bibparskip}%
}
\usepackage{csquotes}
\usepackage[
    backend=biber,
    style=alphabetic, % Or numeric/authoryear
    sorting=nyt,
    doi=true,
    isbn=false,
    url=false,
    ]{biblatex}
\AtEveryBibitem{\clearfield{urlyear}\clearfield{pubstate}} % Suppress access dates and publication status

\theoremstyle{plain}
\newtheorem{theorem}{Theorem}[section]
\newtheorem{introtheorem}{Theorem}

\crefname{introtheorem}{Theorem}{Theorems}
\Crefname{introtheorem}{Theorem}{Theorems}

\newtheorem{question}[theorem]{Question}

\newtheorem{proposition}[theorem]{Proposition}
\newtheorem{corollary}[theorem]{Corollary}

\theoremstyle{definition}

\newtheorem{definition}[theorem]{Definition}
\newtheorem{remark}[theorem]{Remark}

\begin{document}
\title{On explicit formulas for Grassmannian polylogarithms and Gelfand-MacPherson forms}
\author{Juan Diego Rojas}
\address{Department of Mathematics, Northwestern University, Evanston, IL 60208}
\email{reds@u.northwestern.edu}
\begin{abstract}
  We answer a question of Goncharov \Cite{goncharovPolylogarithmsRegulatorsArakelov2004} by showing that his Grassmannian $n$-logarithm can be recovered by integrating the Borel invariant form on the symmetric space $\mathrm{SL}_{n}(\mathbf{C})/\mathrm{SU}(n)$ over a naturally defined ideal simplex. As corollaries, we recover the finiteness of primitive canonical integrals of Brown \cite{brownInvariantDifferentialForms2021} and the evaluation for wheel graphs of Brown and Schnetz \cite{brownWheelClassesLocally2025}. We also compute explicitly the Gelfand-MacPherson forms \Cite{gelfandGeometryGrassmanniansGeneralization1982} for powers of the first Pontryagin class and give a general integral formula for the Gelfand-MacPherson forms associated to each component of the Pontryagin character. 
\end{abstract}
\maketitle
\setcounter{section}{-1}
\section{Introduction}\label{sec:introduction}
The dilogarithm is a special function with a long history. It shows up in several areas of mathematics and physics (for a survey see \Cite{zagierDilogarithmFunction2007}). In this paper, we are particularly interested in its appearance in two different geometric contexts. First, it computes the volume of an ideal tetrahedron in hyperbolic $3$-space, in the form of the Bloch-Wigner dilogarithm. Second, it plays a key role in the construction of a local combinatorial formula for the first Pontryagin class, in the form of the Rogers dilogarithm. The single-valued Grassmannian $n$-logarithm introduced by Goncharov \cite{goncharovChowPolylogarithmsRegulators1995,goncharovPolylogarithmsRegulatorsArakelov2004} is a generalisation of the first instance.  On the other hand, Gelfand and MacPherson \cite{gelfandGeometryGrassmanniansGeneralization1982} introduced the Gelfand-MacPherson forms in an attempt to generalise the second one. The main goal of this paper is to study explicit formulas for both of these generalisations.
\subsection{Goncharov's Grassmannian $n$-logarithm}
In \cite{beilinsonNotesMotivicCohomology1987}, Beilinson, MacPherson, and Schechtman posed as a question the existence of a certain explicit Grassmannian cocycle representing the Beilinson regulator \Cite{beilinsonHigherRegulatorsValues1985} from algebraic $K$-theory to Deligne cohomology (see also \cite{hainHigherLogarithms1990}). In this paper we focus on Goncharov's single-valued Grassmannian $n$-logarithm $\mathcal{L}_{n}^{\mathrm{G}}$ \cite{goncharovChowPolylogarithmsRegulators1995,goncharovGrassmannianTrilogarithms2001,goncharovPolylogarithmsRegulatorsArakelov2004}, which may be extended to a Grassmannian cocycle as conjectured in \Cite{beilinsonNotesMotivicCohomology1987}. Goncharov's Grassmannian $n$-logarithm $\mathcal{L}_{n}^{\mathrm{G}}$ is a function defined on configurations of $2n$ hyperplanes in $\mathbf{CP}^{n-1}$. For $n=2$, $\mathcal{L}_{2}^{\mathrm{G}}$ recovers the Bloch-Wigner dilogarithm of the cross-ratio of four points in $\mathbf{CP}^{1}$. 

The interpretation of the Bloch-Wigner dilogarithm as the volume of an ideal tetrahedron in hyperbolic $3$-space suggests looking at the following geometric generalisation, which Goncharov \Cite[\textsection 7.4]{goncharovPolylogarithmsRegulatorsArakelov2004} attributes to Kontsevich and Nekov\'a\v{r} independently. Let $\mathcal{H}_{n}$ be the symmetric space of positive definite Hermitian forms on $\mathbf{C}^{n}$. Recall that $\mathcal{H}_{n}$ admits a compactification $\overline{\mathcal{H}}_{n}$ (the Satake compactification) realised by nonzero positive-semidefinite forms. Let $h_{1}, \ldots, h_{2n} \subset \mathbf{CP}^{n-1}$ be hyperplanes with no point in common. Each hyperplane $h_{i}$ defines a positive ray of rank-one positive-semidefinite Hermitian forms whose kernels are exactly $h_{i}$. Let $C(h_{1}, \ldots, h_{2n})$ be the projectivisation of the convex hull of the rays corresponding to these hyperplanes. We regard this as a possibly degenerate ideal $(2n-1)$-simplex. Its vertices lie on $\partial \overline{\mathcal{H}}_{n}$, and the interior is contained in $\mathcal{H}_{n}$. The symmetric space $\mathcal{H}_{n}$ admits a naturally defined differential $(2n-1)$-form $\omega_{n}$ called the Borel invariant form. If we identify $\mathcal{H}_{n}$ with $\mathrm{SL}_{n}(\mathbf{C})/\mathrm{SU}(n)$, then $\omega_{n}$ is induced by $\mathrm{tr}((g^{-1}\mathrm{d}g)^{2n-1})$ on $\mathrm{SL}_{n}(\mathbf{C})$. In \Cite{borelStableRealCohomology1974} Borel showed that these differential forms generate the stable cohomology of the general linear group and they are used to define the Borel regulator in algebraic $K$-theory \cite{borelCohomologie$SL_n$Valeurs1977}. In \Cite[\textsection 7.4]{goncharovPolylogarithmsRegulatorsArakelov2004}, Goncharov asks whether the integral of $\omega_{n}$ over $C(h_{1}, \ldots, h_{2n})$ converges and how it compares with $\mathcal{L}_{n}^{\mathrm{G}}(h_{1}, \ldots, h_{2n})$. Our first result is an answer to his question.
\begin{introtheorem}[\Cref{cor:integral-Borel-equals-grassmannian}]\label{thm:introduction-Grassmannian}
  Let $n\geq 2$ and let $h_{1},\ldots,h_{2n}\subset\mathbf{CP}^{n-1}$ be hyperplanes with no point in common. The integral of the Borel invariant form $\omega_{n}$ over $C(h_{1},\ldots,h_{2n})$ converges absolutely and
  \[
    \int_{C(h_{1},\ldots,h_{2n})}\omega_{n}
    = -\frac{2n-1}{4^{n-1}}\binom{2n-2}{n-1}\mathcal{L}_{n}^{\mathrm{G}}(h_{1},\ldots,h_{2n}).
  \]
\end{introtheorem}
Goncharov asks the analogous question \Cite[\textsection 7.3]{goncharovPolylogarithmsRegulatorsArakelov2004} for the geodesic simplex. In \cref{rem:geodesic-vanishing}, we show that this question has a negative answer already at $n=3$. 

\Cref{thm:introduction-Grassmannian} has the following two applications. First, expanding out the trace form $\omega_{n}$ in terms of matrix minors gives the formula obtained by Dan \Cite[Theorem 4.1]{danMultitransgressionRegulators2008} through multitransgression of the Chern character. This answers his question on the identification of his formula with Goncharov's Grassmannian $n$-logarithm $\mathcal{L}_{n}^{\mathrm{G}}$. The further identification with the maximal period of the associated mixed Hodge structure is not treated here.

Second, let $G$ be a connected, bridgeless graph with $2n$ ordered edges and $n$ independent cycles. In \Cite{brownInvariantDifferentialForms2021}, Brown studied the integrals $I_{G}(\omega^{2n-1})$ of the Borel invariant form $\omega_{n}$ over a $(2n-1)$-dimensional simplex $\sigma_{G}$ associated to the graph $G$ and showed that they are finite. These integrals are called primitive canonical integrals, they are a particular instance of the Feynman integrals that appear in quantum field theory. As a corollary of \cref{thm:introduction-Grassmannian}, we recover the finiteness of primitive canonical integrals and compute their values in terms of $\mathcal{L}_{n}^{\mathrm{G}}$. More precisely, the edges $e_{1}, \ldots, e_{2n}$ of $G$ determine hyperplanes $h_{1},\ldots,h_{2n}$ and we show that
\[
  I_{G}(\omega^{2n-1}) = -n\binom{2n}{n}\mathcal{L}_{n}^{\mathrm{G}}(h_{1},\ldots,h_{2n}).
\]
Let $n$ be odd. For the wheel graphs $W_{n}$, we use Goncharov's evaluation of the Grassmannian $n$-logarithm on special configurations \Cite[Theorem 4.5]{goncharovPolylogarithmsRegulatorsArakelov2004} to recover the formula
\[
  I_{W_{n}}(\omega^{2n-1}) = n\binom{2n}{n}\zeta(n)
\]
of Brown and Schnetz \Cite[Theorem 1]{brownWheelClassesLocally2025}. In particular, we obtain a new proof of the non-vanishing of the corresponding wheel classes in the homology of Kontsevich's graph complex. Patzt and Willwacher also prove the non-vanishing of the wheel classes using Goncharov's Grassmannian $n$-logarithm and its evaluation on special configurations \cite[Theorem 7.4 and Remark 7.7]{patztGraphComplexHomology2026}. 

After the first version of the present article appeared on the arXiv, it was pointed out to the author in private communication the resemblance between \cref{thm:introduction-Grassmannian} and \Cite[Theorem 2.2]{portnerGraphIntegralsFeynman2026}. Indeed, Portner's theorem \Cite[Theorem 2.2]{portnerGraphIntegralsFeynman2026} identifies the LHS of \cref{thm:introduction-Grassmannian} (for real configurations) with the Rossi-Willwacher integrals \cite{rossiEtingofsConjectureDrinfeld2014} from deformation quantisation. 

\subsection{Gelfand-MacPherson forms} The appearance of the Rogers dilogarithm in the combinatorial formulas of Gabrielov, Gelfand and Losik for the first Pontryagin class \Cite{gabrielovCombinatorialCalculusCharacteristic1975,gabrielovCombinatorialComputationCharacteristic1975} led Gelfand and MacPherson \cite{gelfandGeometryGrassmanniansGeneralization1982} to the following geometric generalisation. For a polynomial $P$ in Pontryagin classes, they integrate its harmonic representative on the real Grassmannian $\mathrm{G}_{p}^{q}(\mathbf{R}) = \mathrm{Gr}(p, \mathbf{R}^{p+q})$ along the orbits of the Cartan subgroup of $\mathrm{PGL}_{p+q}(\mathbf{R})$. The resulting form $\widehat{P}_{p}^{q}$, which we call Gelfand-MacPherson form, lives on a space of projective configurations $\mathrm{C}_{p}^{q}(\mathbf{R})$. For $P=p_{1}$ on $\mathrm{G}_{2}^{2}(\mathbf{R})$, the resulting one-form $\widehat{P}_{2}^{2}$ is a multiple of the differential of Rogers dilogarithm. We compute explicitly the forms $\widehat{P}_{2}^{2k}$ for $P=p_{1}^{k}$, the $k$th power of the first Pontryagin class.
\begin{introtheorem}[\Cref{thm:Gelfand-MacPherson-polylogarithmic-form-rank-2}]\label{thm:introduction-GM-rank-two}
  Let $P=p_{1}^{k}$ be the $k$th power of the first Pontryagin class. In coordinates $x = (x_{1},\ldots,x_{2k+1}) \in \mathrm{Conf}_{2k+1}(\mathbf{R})/\mathrm{Aff}_{1}(\mathbf{R}) \cong \mathrm{C}_{2}^{2k}(\mathbf{R})$, we have that $\widehat{P}_{2}^{2k} \in \Omega^{2k-1}(\mathrm{C}_{2}^{2k}(\mathbf{R}))$ is given by
  \[
     \widehat{P}_{2}^{2k} = -\frac{4}{(2\pi)^{2k}}\sum_{i<j}(-1)^{i+j}\log|x_{i}-x_{j}|\bigwedge_{\ell\neq i,j}\mathrm{d}\log \left|\frac{x_{\ell} - x_{i}}{x_{\ell} - x_{j}}\right|.
  \]
\end{introtheorem}
For $k=1$, the formula recovers the differential of Rogers dilogarithm. For larger $k$, it verifies a conjecture of Pirio, which he had verified up to $k \leq 6$ \cite[Proposition 5.28]{pirioN+3websRationalCurves2023}. It was pointed out to the author that Pirio settled the formula on \Cite{pirioEulersAbelianRelation2026} shortly before the first version of the present paper.

We also show how Brown's algebra of multiple polylogarithms \cite{brownMultipleZetaValues2009} can be used to transgress the form $\widehat{P}_{2}^{2k}$ to a linear combination of multiple polylogarithms in $\mathcal{M}_{0,4k}(\mathbf{R})$. We ask whether it is always the case that the transgressed function can be written using only Rogers dilogarithms and logarithms.

For the $k$th component of the Pontryagin character $\mathrm{ph}_{k}$, Damiano's vanishing theorem \Cite{damianoVanishingCharacteristicForms1982} singles out the self-dual Grassmannian $\mathrm{G}_{2k}^{2k}(\mathbf{R})$. Among the Grassmannians $\mathrm{G}_{p}^{q}(\mathbf{R})$ with $p+q=4k$, it is the only possible source of a nonzero Gelfand-MacPherson one-form associated to $\mathrm{ph}_{k}$. We give an integral formula for $\widehat{P}_{2k}^{2k}$ in this case.
\begin{introtheorem}[\Cref{thm:self-dual-GM-simplex-integrand}]\label{thm:introduction-GM-self-dual}
  Let $k\geq 1$ and $P=\mathrm{ph}_{k}$. Represent a generic configuration in $\mathrm{C}_{2k}^{2k}(\mathbf{R})$ by
  \[
    S(x) = (I_{2k}\mid\mathbf{1}\mid z_{1}\mid\cdots\mid z_{2k-1}) = (v_{1},\ldots,v_{4k}),
    \qquad z_{b} = (x_{1b},\ldots,x_{2k-1,b},1)^{t}.
  \]
  Here $\mathbf{1}$ is the column whose $2k$ entries are all $1$. Set
  \[
    Q(r,x) = \det\left(\sum_{i=1}^{4k}r_{i}v_{i}v_{i}^{t}\right)
    = \sum_{|I|=2k}r_{I}\Delta_{I}(S(x))^{2}, \qquad r_{I}=\prod_{i\in I}r_{i},
  \]
  where $\Delta_{I}$ denotes the Plücker coordinate indexed by $I$. Then
  \[
    \widehat{P}_{2k}^{2k}
    = \frac{(-1)^{k-1}2^{4k-1}}{(2k-1)!(2\pi)^{2k}}
    \sum_{a,b=1}^{2k-1}\left(\int_{\Delta^{4k-1}}\frac{N_{ab}(r,x)}{Q(r,x)^{k+1}}\,\mathrm{d}r\right)\mathrm{d}x_{ab},
  \]
  where $N_{ab}(r,x)$ are explicitly determined polynomials defined in \eqref{eqn:self-dual-simplex-numerator}. Here $\mathrm{d}r=\mathrm{d}r_{1}\cdots\mathrm{d}r_{4k-1}$ and $r_{4k}=1-\sum_{i<4k}r_{i}$.
\end{introtheorem}
\subsection{The Genocchi-Hermite interpolation formula}
A key idea that appears throughout the paper is to systematically use the Genocchi-Hermite interpolation formula \Cite{genocchiRelationEntreDifference1869,hermiteFormuleDinterpolationLagrange1877} in simplicial Chern-Weil theory. Applied to the logarithm, it evaluates a simplex integral as a logarithmic cocycle. We expect the following identity to be useful beyond the cases treated in this paper.
\begin{introtheorem}[Logarithmic transgression; \Cref{thm:barycentric-transgression}]\label{thm:introduction-logarithmic-transgression}
  Let $k\geq 1$, let $X$ be a closed oriented $(k-1)$-manifold, and let $a_{0},\ldots,a_{k}\colon X\to\mathbf{R}_{>0}$ be smooth positive functions. For $t\in\Delta^{k}=\{(t_{0},\ldots,t_{k})\in\mathbf{R}_{\geq 0}^{k+1}\colon\sum t_{i}=1\}$, set $A_{t}=\sum t_{i}a_{i}$ and $b_{i}=a_{i}-a_{0}$. Then
  \[
    \int_{\Delta^{k}}\left(\int_{X}\frac{b_{1}}{A_{t}}\mathrm{d}\left(\frac{b_{2}}{A_{t}}\right)\wedge\cdots\wedge\mathrm{d}\left(\frac{b_{k}}{A_{t}}\right)\right)\mathrm{d}t
    = \frac{1}{k!}\int_{X}\log\frac{a_{1}}{a_{0}}\,\mathrm{d}\log\frac{a_{2}}{a_{0}}\wedge\cdots\wedge\mathrm{d}\log\frac{a_{k}}{a_{0}},
  \]
  where $\mathrm{d}$ differentiates on $X$ and $\mathrm{d}t=\mathrm{d}t_{1}\cdots\mathrm{d}t_{k}$. As a function of $(a_{0},\ldots,a_{k})$, this integral is a simplicial $k$-cocycle: its alternating sum over the tuples obtained by omitting one of $k+2$ inputs is zero.
\end{introtheorem}
\subsection*{Conventions on Grassmannian polylogarithms} In this paper we follow the terminology of referring to a function as polylogarithmic if it may be written using linear combinations of multiple polylogarithms $\mathrm{Li}_{n_{1}, \ldots, n_{k}}$ with rational functions as arguments. Goncharov's universality conjecture predicts that periods of mixed Tate motives are polylogarithmic.

There are several different constructions of ``Grassmannian polylogarithms'' in the literature. There is a geometric construction due to Hanamura and MacPherson of a Grassmannian $n$-logarithm \Cite{hanamuraGeometricConstructionPolylogarithms1993,hanamuraGeometricConstructionPolylogarithms1996}. The author is unaware whether those functions are known to be polylogarithmic nor how they relate to the rest of constructions. There is the real single-valued Grassmannian $n$-logarithm of Goncharov $\mathcal{L}_{n}^{\mathrm{G}}$ that we study in this paper. The author believes it is open whether these are polylogarithmic beyond weight $3$. There is the complex-analytic multi-valued function $\mathrm{Gr}_{n}$ constructed by Goncharov in \Cite{goncharovSimpleConstructionGrassmannian2013}. These were shown to be polylogarithmic by Charlton, Gangl, and Radchenko \Cite{charltonAomotoGrassmannianPolylogarithms2025a}. It is unknown what the exact relationship between $\mathcal{L}_{n}^{\mathrm{G}}$ and $\mathrm{Gr}_{n}$ is \Cite[\textsection 1.2]{goncharovSimpleConstructionGrassmannian2013}. There is the construction of Gelfand and MacPherson where they associate differential forms to Pontryagin polynomials, as we recall in \cref{sec:GM-polylogarithms}. It is unknown whether these forms are polylogarithmic beyond the case of \cref{thm:introduction-GM-rank-two} where we explicitly compute them. In \cite[\textsection 1.2]{goncharovSimpleConstructionGrassmannian2013}, Goncharov mentions that Gelfand and MacPherson define a real-valued Grassmannian $2n$-logarithm by constructing its differential. For $n\geq2$, this description requires justification: the construction produces a one-form whose closedness lies beyond the known vanishing theorems, as explained in \cref{subsec:closedness}. Until this is established, the construction does not define the claimed ``Gelfand-MacPherson Grassmannian polylogarithm'', even locally.
\subsection*{Notation} We use boldface for the standard objects $\mathbf{N}$, $\mathbf{Z}$, $\mathbf{Q}$, $\mathbf{R}$, and $\mathbf{C}$. All varieties are defined over $\mathbf{R}$ or $\mathbf{C}$. We index Grassmannians by dimension and codimension, i.e.\ $\mathrm{G}_{p}^{q} = \mathrm{Gr}(p,p+q)$. For projective spaces, we use $\mathbf{RP}^{n}$ (resp.\ $\mathbf{CP}^{n}$) when using its structure as a smooth (resp.\ complex) manifold and use $\mathbf{P}^{n}$ for algebraic constructions. 
\subsection*{Acknowledgements} This paper is part of the author's Ph.D.\ thesis at Northwestern University supervised by Daniil Rudenko and Boris Tsygan. He thanks them for all the encouragement and impeccable supervision throughout his graduate studies. The author would like to thank Daniil Rudenko for recommending the paper \Cite{gelfandGeometryGrassmanniansGeneralization1982} which started this whole project and for believing in him when he decided to write his thesis on the subject. He would also like to thank Ben Antieau, Paul Bressler and Gus Schrader for showing interest in the author's work and for many discussions. The author thanks Paul Bressler for several comments on a previous draft. Finally, all the key ideas in this paper had their genesis while the author was in residence for a season in Buenos Aires and while hosted at Universidad de los Andes in Bogotá. The project would have reached a dead end if it had not been for the excellent working conditions and vibrant life in both cities. 
\subsection*{AI disclosure}
This project originated before AI tools were developed enough to meaningfully contribute. In the last stage, the author used a personal ChatGPT Plus subscription to proof-read the final draft and help with citation management, literature search, and LaTeX. Proof-reading solved a previous issue with the constants in \cref{prop:value-grassmannian-special-configuration}, since the published version of \cite{goncharovPolylogarithmsRegulatorsArakelov2004} contains a typo that is not present in \cite{levinNoteRHodgeTateSheaves2001}. 
\section{Simplex integration and transgression: an application of the Genocchi-Hermite interpolation formula}\label{sec:simplex-integration}
In this section we study how the Genocchi-Hermite interpolation formula can be applied to simplicial Chern-Weil theory \cite{dupontSimplicialRhamCohomology1976} and the transgression identities we can derive from it. Let
\[
  \Delta^{m} = \left\{ (t_{0}, \ldots, t_{m}) \in \mathbf{R}^{m+1}_{\geq 0} \colon \sum t_{i} = 1 \right\}
\]
be the standard simplex, oriented by $\mathrm{d}t_{1}\wedge\cdots\wedge\mathrm{d}t_{m}$. Let $Y_{\bullet}$ be a simplicial manifold with face maps $\partial_{i}\colon Y_{m}\to Y_{m-1}$. Write $\iota_{i}\colon\Delta^{m-1}\to\Delta^{m}$ for the inclusion of the face $\{t_{i}=0\}$ and let $\delta=\sum_{i}(-1)^{i}\partial_{i}^{*}$ denote the simplicial differential.
\begin{proposition}[Simplex integration]\label{prop:simplex-integration}
  Let $N \in \mathbf{N}$. For each $m \in \mathbf{N}$, let $\eta_{m}\in\Omega^{N}(\Delta^{m}\times Y_{m})$ be smooth closed forms satisfying
  \[
    (\iota_{i}\times\mathrm{id})^{*}\eta_{m} = (\mathrm{id}\times\partial_{i})^{*}\eta_{m-1}.
  \]
  Then the forms
  \[
    T_{m} := \int_{\Delta^{m}}\eta_{m}\in\Omega^{N-m}(Y_{m})
  \]
  satisfy $\mathrm{d}T_{0}=0$ and
  \begin{equation}\label{eqn:simplex-integration-descent}
    \delta T_{m-1}+(-1)^{m}\mathrm{d}T_{m}=0 \qquad (m\geq 1).
  \end{equation}
  In particular, $T_{N}$ is a cocycle, i.e.\ $\delta T_{N}=0$.
\end{proposition}
\begin{proof}
  Fibrewise Stokes' theorem gives
  \[
    \int_{\Delta^{m}}\mathrm{d}\eta_{m}
    = \int_{\partial\Delta^{m}}\eta_{m}+(-1)^{m}\mathrm{d}\int_{\Delta^{m}}\eta_{m}.
  \]
  The left-hand side is zero. Compatibility on faces identifies the boundary integral with $\delta T_{m-1}$. For $m=N+1$, the form $T_{N+1}$ has negative degree and is zero, so $\delta T_{N}=0$.
\end{proof}

\begin{remark}\label{rem:simplicial-Chern-Weil}
  For characteristic classes, the first case of \cref{prop:simplex-integration} is the Chern-Simons transgression. Let $E\to B$ be a vector bundle and let $P$ be an invariant polynomial of degree $r$. Given connections $\nabla_{0},\ldots,\nabla_{m}$ on $E$, form the total connection
\[
  \widetilde\nabla = \mathrm{d}_{\Delta^{m}}+\sum_{i=0}^{m}t_{i}\nabla_{i}
\]
on the pullback of $E$ to $\Delta^{m}\times B$ by the projection. Set
\[
  T_{m}(\nabla_{0},\ldots,\nabla_{m})
  := \int_{\Delta^{m}}P(F_{\widetilde\nabla})\in\Omega^{2r-m}(B),
\]
where $F_{\widetilde\nabla}$ denotes the curvature. The Bianchi identity implies that $P(F_{\widetilde\nabla})$ is closed. If we restrict to a face, we just omit the connection corresponding to that face. Therefore, \eqref{eqn:simplex-integration-descent} holds. For two connections, we obtain
\[
  \mathrm{d}T_{1}(\nabla_{0},\nabla_{1})
  = P(F_{\nabla_{1}})-P(F_{\nabla_{0}}).
\]
The Chern-Simons form $T_{1}$ gives an explicit primitive for the difference of the two characteristic forms. For three connections, we get that
\[
  \mathrm{d}T_{2}(\nabla_{0},\nabla_{1},\nabla_{2})
  = T_{1}(\nabla_{0},\nabla_{2})-T_{1}(\nabla_{0},\nabla_{1})-T_{1}(\nabla_{1},\nabla_{2}).
\]
Thus the transgression from $\nabla_{0}$ to $\nabla_{2}$ differs from the sum of the transgressions from $\nabla_{0}$ to $\nabla_{1}$ and from $\nabla_{1}$ to $\nabla_{2}$ by an exact form; its primitive is the integral over the $2$-simplex. Applying this construction to compatible connections on a simplicial bundle gives the characteristic forms in Dupont's simplicial Chern-Weil theory \Cite{dupontSimplicialRhamCohomology1976}.
\end{remark}
We now focus on computing the integrals in \cref{prop:simplex-integration} in two explicit cases. We first recall the notion of divided differences. Let $f$ be a function defined on an interval $[a,b]$ containing distinct points $x_{0}, \ldots, x_{k}$. Its \textsf{$k$th divided difference} on the points $x_{0}, \ldots, x_{k}$ is defined by the rule
\begin{equation}\label{eqn:divided-difference-definition}
  [x_{0}, \ldots, x_{k}]f := \sum_{i=0}^{k}\frac{f(x_{i})}{\prod_{j\neq i} (x_{i} - x_{j})}.
\end{equation}
Equivalently, divided differences are given recursively by setting $[x_{0}]f := f(x_{0})$ and for $k\geq 1$
\begin{equation}\label{eqn:divided-difference-recursion}
  [x_{0}, \ldots, x_{k}]f := \frac{[x_{1}, \ldots, x_{k}]f-[x_{0},\ldots,x_{k-1}]f}{x_{k}-x_{0}}.
\end{equation}
The Genocchi-Hermite interpolation formula \Cite{genocchiRelationEntreDifference1869,hermiteFormuleDinterpolationLagrange1877} writes the $k$th divided difference of $f$ as the integral over the $k$-simplex of the $k$th derivative of $f$. We include a proof for completeness.
\begin{theorem}[Genocchi-Hermite interpolation formula]\label{thm:hermite-genocchi-formula}
  If $f\colon[a,b] \to \mathbf{R}$ is of class $C^{k}$ then
\begin{equation}\label{eqn:Hermite-Genocchi}
  [x_{0}, \ldots, x_{k}]f = \int_{\Delta^{k}} f^{(k)}\left(\sum_{i=0}^{k} t_{i}x_{i} \right) \mathrm{d}t.
\end{equation}
\end{theorem}
\begin{proof}
  We show that
  \[
    I_{k}(x_{0},\ldots,x_{k}) := \int_{\Delta^{k}} f^{(k)}\left(\sum_{i=0}^{k} t_{i}x_{i} \right) \mathrm{d}t
  \]
  satisfies the recursion \eqref{eqn:divided-difference-recursion}. Write
  \[
    X = x_{0} + \sum_{i=1}^{k-1}t_{i}(x_{i}-x_{0}),\quad T = \sum_{i=1}^{k-1}t_{i}
  \]
  and observe that $0 \leq t_{k} \leq 1-T$ and that $\sum_{i=0}^{k}t_{i}x_{i} = X + t_{k}(x_{k}-x_{0})$. Therefore,
  \[
    I_{k}(x_{0},\ldots,x_{k}) = \int_{\Delta^{k-1}}\left(\int_{0}^{1-T}f^{(k)}(X + t_{k}(x_{k}-x_{0}))\mathrm{d}t_{k}\right) \mathrm{d}t_{1}\cdots \mathrm{d}t_{k-1}.
  \]
  By the fundamental theorem of calculus, we have that
  \[
  (x_{k}-x_{0})\int_{0}^{1-T}f^{(k)}(X + t_{k}(x_{k}-x_{0}))\mathrm{d}t_{k} = f^{(k-1)}(X + (1-T)(x_{k}-x_{0})) - f^{(k-1)}(X).
  \]
  Observe that
  \[
    X + (1-T)(x_{k}-x_{0}) = \sum_{i=1}^{k-1}t_{i}x_{i} + (1-T)x_{k}, \quad X = (1-T)x_{0} + \sum_{i=1}^{k-1}t_{i}x_{i}.
  \]
  Therefore, integrating over $\Delta^{k-1}$ yields
  \[
    (x_{k}-x_{0})I_{k}(x_{0}, \ldots, x_{k}) = I_{k-1}(x_{1}, \ldots, x_{k}) - I_{k-1}(x_{0}, \ldots, x_{k-1})
  \]
  as desired.
\end{proof}
The integral also defines divided differences when the points coincide. We use this convention below. For $f(x)=\log x$, $k\geq 1$ and positive $x_{0},\ldots,x_{k}$, \eqref{eqn:Hermite-Genocchi} gives
\begin{equation}\label{eqn:Hermite-Genocchi-logarithm}
  \int_{\Delta^{k}}\frac{\mathrm{d}t}{(\sum_{i=0}^{k}t_{i}x_{i})^{k}} = \frac{(-1)^{k+1}}{(k-1)!}[x_{0}, \ldots, x_{k}]\log x.
\end{equation}
We now apply this identity to a family of positive functions on a closed oriented $(k-1)$-manifold $X$. Let $a\colon X\to\mathbf{R}_{>0}$ be a smooth positive function and let $b_{1},\ldots,b_{k}\in C^{\infty}(X,\mathbf{R})$. We regard $a$ as a point in the space of positive functions and each $b_{i}$ as a tangent vector at $a$, representing the tangent vector $a+\varepsilon b_{i}$. Consider
\[
  \omega_{a}(b_{1},\ldots,b_{k})
  := \int_{X}\frac{b_{1}}{a}\,\mathrm{d}\!\left(\frac{b_{2}}{a}\right)\wedge\cdots\wedge\mathrm{d}\!\left(\frac{b_{k}}{a}\right).
\]
If the reader wishes not to think about calculus in infinite dimensions, they might as well read the statement after pulling back to a finite-dimensional smooth family of positive functions. In the next theorem, we will do just that, where the family is the simplex $A_{t}=\sum t_{i}a_{i}$ and we set $a=A_{t}$ and $b_{i}=\partial A_{t}/\partial t_{i}=a_{i}-a_{0}$ for $1\leq i\leq k$.
\begin{theorem}[Logarithmic transgression]\label{thm:barycentric-transgression}
  Let $k\geq 1$. Let $X$ be a closed oriented manifold of dimension $k-1$ and let $a_{0}, \ldots, a_{k} \colon X \to \mathbf{R}_{> 0}$ be smooth positive functions. Set $A_{t} = \sum t_{i}a_{i}$ with $t = (t_{0}, \ldots, t_{k}) \in \Delta^{k}$. Set $b_{i} = a_{i} - a_{0}$, then the form $\omega$ above is closed and
  \begin{equation}\label{eqn:barycentric-transgression}
    \int_{\Delta^{k}}\left(\int_{X} \frac{b_{1}}{A_{t}}\mathrm{d}\left(\frac{b_{2}}{A_{t}}\right) \wedge \cdots \wedge \mathrm{d}\left(\frac{b_{k}}{A_{t}}\right) \right) \mathrm{d}t = \frac{1}{k!}\int_{X} \log \frac{a_{1}}{a_{0}} \mathrm{d}\log \frac{a_{2}}{a_{0}} \wedge \cdots \wedge \mathrm{d}\log \frac{a_{k}}{a_{0}}.
  \end{equation}
  Equivalently, the left-hand side is $\int_{\Delta^{k}}A_{t}^{*}\omega$. As a function of the tuple $(a_{0},\ldots,a_{k})$, either side is a simplicial $k$-cocycle, with face maps given by forgetting an input.
\end{theorem}
\begin{proof}
We first check that $\omega$ is alternating and closed. The expression
\[
  (u_{1},\ldots,u_{k})\longmapsto\int_{X}u_{1}\,\mathrm{d}_{X}u_{2}\wedge\cdots\wedge\mathrm{d}_{X}u_{k}
\]
is alternating in the last $k-1$ entries. For $k\geq 2$, Stokes' theorem gives
\[
  \int_{X}u_{1}\,\mathrm{d}_{X}u_{2}\wedge\cdots\wedge\mathrm{d}_{X}u_{k}
  =-\int_{X}u_{2}\,\mathrm{d}_{X}u_{1}\wedge\mathrm{d}_{X}u_{3}\wedge\cdots\wedge\mathrm{d}_{X}u_{k},
\]
so it is alternating in all entries. In the coordinate $u=\log a$, a tangent vector $b$ becomes $b/a$. Thus $\omega$ has constant coefficients in this coordinate and is closed. The same conclusion is immediate for $k=1$.

We next evaluate its integral over the simplex. For $k=1$, the identity is the fundamental theorem of calculus applied to $\log A_{t}$. Assume that $k\geq 2$. For smooth functions $c_{1}, \ldots, c_{k}$ on $X$, define the radial $(k-1)$-form
\[
 \rho(c) = \sum_{i=1}^{k} (-1)^{i+1}c_{i} \, \mathrm{d}c_{1} \wedge \cdots \wedge \widehat{\mathrm{d}c_{i}} \wedge \cdots \wedge \mathrm{d}c_{k}. 
\]
Stokes' gives $\int_{X} \rho(c) = k \int_{X} c_{1} \mathrm{d}c_{2} \wedge \cdots \wedge \mathrm{d}c_{k}$. Set $c_{i} = \frac{b_{i}}{A_{t}}$, so that the left-hand side of \eqref{eqn:barycentric-transgression} equals
\begin{equation}\label{eqn:lhs-hermite-manifold}
  \frac{1}{k}\int_{X} \rho(b) \int_{\Delta^{k}} A_{t}^{-k} \mathrm{d}t.
\end{equation}
By the Genocchi-Hermite interpolation formula for the logarithm \eqref{eqn:Hermite-Genocchi-logarithm}, we may write \eqref{eqn:lhs-hermite-manifold} as
\begin{equation}\label{eqn:post-hermite}
  \frac{(-1)^{k+1}}{k!}\int_{X}[a_{0}, \ldots, a_{k}]\log \rho(b).
\end{equation}
To evaluate \eqref{eqn:post-hermite}, set $x_{i}=a_{i}/a_{0}$ and consider the map $x\colon X\to\mathbf{R}^{k}_{>0}$. On this positive orthant, write $y_{i}=x_{i}-1$ and $F=[1,x_{1},\ldots,x_{k}]\log$. Observe that
\[
  [a_{0},\ldots,a_{k}]\log\,\rho(b)=x^{*}\bigl(F\rho(y)\bigr).
\]
We will prove on the positive orthant that
\begin{equation}\label{eqn:claim-modulo}
  (-1)^{k+1}F\rho(y)
  \equiv \log x_{1}\,\mathrm{d}\log x_{2}\wedge\cdots\wedge\mathrm{d}\log x_{k}
  \mod \mathrm{d}\Omega^{k-2}(\mathbf{R}^{k}_{>0}).
\end{equation}
Let $E=\sum y_{i}\partial_{y_{i}}$ denote the Euler vector field. Then
\[
  \mathrm{d}(F\rho(y))=(EF+kF)\,\mathrm{d}y_{1}\wedge\cdots\wedge\mathrm{d}y_{k}.
\]
To compute the coefficient, write $B_{t}=t_{0}+\sum_{i=1}^{k}t_{i}x_{i}$. Formula \eqref{eqn:Hermite-Genocchi-logarithm} and $EB_{t}=B_{t}-1$ give
\[
  (E+k)F=(-1)^{k+1}k!\int_{\Delta^{k}}B_{t}^{-k-1}\,\mathrm{d}t
  =\frac{(-1)^{k+1}}{x_{1}\cdots x_{k}}.
\]
For the last equality, use the change of variables $s_{0}=t_{0}/B_{t}$ and $s_{i}=t_{i}x_{i}/B_{t}$. This change of variables has Jacobian $(x_{1}\cdots x_{k})B_{t}^{-k-1}$. The $k!$ comes from the reciprocate of the volume of the simplex. Therefore, both forms in \eqref{eqn:claim-modulo} have differential $\mathrm{d}\log x_{1}\wedge\cdots\wedge\mathrm{d}\log x_{k}$. Their difference is closed and has positive degree. Since the positive orthant is contractible, this difference is exact. Pulling back to $X$ and integrating proves \eqref{eqn:barycentric-transgression}. The cocycle statement is immediate from the simplicial structure.
\end{proof}
\begin{theorem}[Quadratic simplex integral]\label{thm:quadratic-barycentric-transgression}
  Let $k\geq 1$ and let $x_{0}, \ldots, x_{2k} \in \mathbf{R}$ be pairwise distinct. Then
  \[
    \int_{\Delta^{2k}}\frac{\mathrm{d}t}{\left(\sum_{i<j}t_{i}t_{j}(x_{i}-x_{j})^{2}\right)^{k}} = \frac{4k(-1)^{k}}{(2k)!}\sum_{i < j}\frac{(x_{i}-x_{j})^{2k-2}\log|x_{i}-x_{j}|}{\prod_{\ell \neq i,j}(x_{i}-x_{\ell})(x_{j}-x_{\ell})}  
  \] 
\end{theorem}
\begin{proof}
Set $V(t) := \sum_{i<j}t_{i}t_{j}(x_{i}-x_{j})^{2}$. Observe that for any $y$,
\[
  \sum_{i}t_{i}(y-x_{i})^{2} = \left(y - \sum_{i}t_{i}x_{i}\right)^{2} + V(t).
\]
Using the beta integral
\[
  \int_{\mathbf{R}} (1+u^{2})^{-k-\frac{1}{2}}\mathrm{d}u = \frac{4^{k}(k!)^{2}}{k(2k)!},
\]
one may see that
\[
  V(t)^{-k} = \frac{k(2k)!}{4^{k}(k!)^{2}}\int_{\mathbf{R}}\left(\sum_{i}t_{i}(y-x_{i})^{2}\right)^{-k- \frac{1}{2}} \mathrm{d}y.
\]
We now apply the Genocchi-Hermite interpolation formula to the function $f(s)=s^{k- \frac{1}{2}}$ at the $2k+1$ points $(y-x_{0})^{2}, \ldots, (y-x_{2k})^{2}$. Observe that
\[
  f^{(2k)}(s) = (-1)^{k}\frac{(2k)!^{2}}{2^{4k}(k!)^{2}}s^{-k- \frac{1}{2}},
\]
so that
\[
  \int_{\Delta^{2k}} \frac{\mathrm{d}t}{V(t)^{k}} = \frac{(-1)^{k}4^{k}k}{(2k)!}\int_{\mathbf{R}}F(y)\mathrm{d}y
\]
where
\[
  F(y) = [(y-x_{0})^{2},\ldots,(y-x_{2k})^{2}]s^{k- \frac{1}{2}}.
\]
From the explicit formula \eqref{eqn:divided-difference-definition} we see that
\[
  F(y) = \sum_{i}\frac{|y-x_{i}|^{2k-1}}{\prod_{j\neq i}(x_{j}-x_{i})\prod_{j\neq i}(2y-x_{i}-x_{j})}.
\]
Since $2k-1$ is odd, writing $\varepsilon_{i}(y) = \mathrm{sgn}(y-x_{i})$ gives
\[
  |y-x_{i}|^{2k-1} = \varepsilon_{i}(y)(y-x_{i})^{2k-1}.
\]
Partial fractions in $y$ gives
\[
  F(y)=- \frac{1}{2^{2k-1}}\sum_{i<j}\frac{(x_{i}-x_{j})^{2k-2}}{\prod_{\ell\neq i,j}(x_{i}-x_{\ell})(x_{j}-x_{\ell})}\frac{\varepsilon_{i}(y) + \varepsilon_{j}(y)}{2y-x_{i}-x_{j}}.
\]
The remaining calculation is to carry out the integral $\int_{\mathbf{R}} F(y)\mathrm{d}y$. Integrating separately on the intervals determined by $x_{i}$ and $x_{j}$ gives
\[ 
  \int_{-R}^{R} \frac{\varepsilon_{i}(y) + \varepsilon_{j}(y)}{2y - x_{i}-x_{j}}\mathrm{d}y = 2\log(2R) - 2\log|x_{i}-x_{j}| + o(1).
\]
Summing up, the terms $2\log(2R)$ cancel. Since $F(y) = O(|y|^{-2k-1})$, taking $R \to \infty$ gives
\[
  \int_{\mathbf{R}} F(y)\mathrm{d}y = 2^{2-2k}\sum_{i < j} \frac{(x_{i}-x_{j})^{2k-2}\log|x_{i}-x_{j}|}{\prod_{\ell \neq i,j}(x_{i}-x_{\ell})(x_{j}-x_{\ell})}
\]
which implies the desired result. 
\end{proof}
\section{Goncharov's Grassmannian $n$-logarithm}\label{sec:grassmannian-polylogarithms}
The main goal of this section is to answer a question of Goncharov in  \Cite[\textsection 7.4]{goncharovPolylogarithmsRegulatorsArakelov2004} regarding a different approach to the single-valued Grassmannian $n$-logarithm $\mathcal{L}_{n}^{\mathrm{G}}$. Let $n \geq 2$.
% Let $X$ be a complex variety.  For nonzero rational functions $f_{1}, \ldots, f_{m}$ on a complex variety, consider the $\mathbf{R}(n-1)$-valued $(n-1)$-form (see \Cite[Sec. 1.2.]{goncharovGrassmannianTrilogarithms2001})
% \begin{equation}\label{eqn:singular-r-form}
%   r_{m}(f_{1}, \ldots, f_{m}) := - \mathrm{Alt}_{m}\sum_{\substack{ j \geq 0 \\ 2j+1 \leq m }} \frac{{m \choose 2j+1}}{m!}\log|f_{1}| \bigwedge_{a=2}^{2j+1}\mathrm{d}\log|f_{a}| \wedge \bigwedge_{a = 2j+2}^{m} \sqrt{-1}\mathrm{d}\mathrm{arg}(f_{a}).
% \end{equation}
\subsection{Setting}
\begin{definition}\label{def:grassmannian-n-logarithm}\Cite[Prop. 6.2.]{goncharovGrassmannianTrilogarithms2001}
Let $h_{1}, \ldots, h_{2n} \subset \mathbf{CP}^{n-1}$ be a configuration of hyperplanes. Choose rational functions $f_{i}$ on $\mathbf{CP}^{n-1}$ such that $\mathrm{div}(f_{i}) = h_{i+1} - h_{1}$ for $1 \leq i \leq 2n-1$. The single-valued \textsf{Grassmannian} $n$-\textsf{logarithm} is given by the formula
\begin{equation}\label{eqn:grassmannian-n-logarithm}
  \mathcal{L}^{\mathrm{G}}_{n}(h_{1}, \ldots, h_{2n}) := - \frac{(-4)^{n-1}( (n-1)! )^{2}}{(2\pi i)^{n-1}(2n-2)!}\int_{\mathbf{CP}^{n-1}} \log|f_{1}|\mathrm{d}\log|f_{2}| \wedge \cdots \wedge \mathrm{d}\log|f_{2n-1}|.
\end{equation}
\end{definition}
The Grassmannian $n$-logarithm satisfies skew-symmetry and a pair of functional equations.
\begin{theorem}[{\Cite[Theorem 4.3]{goncharovPolylogarithmsRegulatorsArakelov2004}}]\label{thm:properties-grassmannian-polylogarithm}
\leavevmode
  \begin{enumerate}
    \item The function $\mathcal{L}_{n}^{\mathrm{G}}$ is skew-symmetric, i.e.\ for any $2n$ hyperplanes $h_{1}, \ldots, h_{2n}$ and any $\sigma \in \Sigma_{2n}$ we have that
    \[
      \mathcal{L}^{\mathrm{G}}_{n}(h_{1}, \ldots, h_{2n}) = (-1)^{|\sigma|}\mathcal{L}_{n}^{\mathrm{G}}(h_{\sigma(1)}, \ldots, h_{\sigma(2n)}).
    \]
    \item For any $2n+1$ hyperplanes $h_{1}, \ldots, h_{2n+1}$ in $\mathbf{CP}^{n-1}$ we have that
    \[
      \sum_{i=1}^{2n+1} (-1)^{i} \mathcal{L}_{n}^{\mathrm{G}}(h_{1}, \ldots, \widehat{h}_{i}, \ldots, h_{2n+1}) = 0.
    \]
    \item For any $2n+1$ hyperplanes in $\mathbf{CP}^{n}$ we have that
    \[
      \sum_{\substack{i=1 \\ i \text{ admissible}}}^{2n+1} (-1)^{i}\mathcal{L}_{n}^{\mathrm{G}}(h_{i} \cap h_{1}, \ldots, \widehat{h_{i} \cap h_{i}}, \ldots, h_{i} \cap h_{2n+1}) = 0,
    \]
    where $i$ admissible means that $\bigcap_{j\neq i}h_{j} = \emptyset$ and there is no $j \neq i$ with $h_{j} = h_{i}$.\footnote{The admissibility condition was omitted from Goncharov's published statement. This omission, confirmed by Goncharov, is explained in \cite[Remark 7.2]{patztGraphComplexHomology2026}.}
  \end{enumerate}
\end{theorem}
\begin{remark}\label{rem:grassmannian-bloch-wigner}
  For $n=2$, the Grassmannian dilogarithm is the \textsf{Bloch-Wigner dilogarithm}. Recall that
  \[
    D(z) := \operatorname{Im}\mathrm{Li}_{2}(z) + \log|z|\arg(1-z).
  \]
  It extends to a single-valued real-analytic function on $\mathbf{C}\setminus\{0,1\}$ and to a continuous function on $\mathbf{CP}^{1}$ by setting $D(0)=D(1)=D(\infty)=0$. For four distinct points $z_{1},\ldots,z_{4}\in\mathbf{CP}^{1}$, put
  \[
    \mathrm{cr}(z_{1},z_{2},z_{3},z_{4}) := \frac{(z_{1}-z_{3})(z_{2}-z_{4})}{(z_{1}-z_{4})(z_{2}-z_{3})},
  \]
  for the cross-ratio. It follows that
  \[
    \mathcal{L}_{2}^{\mathrm{G}}(z_{1},z_{2},z_{3},z_{4}) = i\,D\bigl(\mathrm{cr}(z_{1},z_{2},z_{3},z_{4})\bigr),
  \]
  as in \Cite[\textsection 4.5, Eq. (21)]{goncharovGrassmannianTrilogarithms2001}. In particular, $\mathcal{L}_{2}^{\mathrm{G}}(\infty,0,1,z)=i\,D(z)$. The differential is
  \[
    \mathrm{d}D(z) = \log|z|\,\mathrm{d}\arg(1-z) - \log|1-z|\,\mathrm{d}\arg z.
  \]
  Moreover, $D(\overline{z})=-D(z)$, so $D$ vanishes on the real projective line. It satisfies the five-term relation 
  \begin{equation}\label{eqn:five-term-relation-cross-ratios}
    \sum_{j=0}^{4}(-1)^{j}D(\mathrm{cr}(z_{0}, \ldots, \widehat{z_{j}}, \ldots, z_{4})) = 0
  \end{equation}
  which is equivalent to both of the functional equations in \cref{thm:properties-grassmannian-polylogarithm}. Geometrically, $D(z)$ is the oriented volume of an ideal hyperbolic tetrahedron with cross-ratio $z$ \Cite[\textsection 1, Eq. (12) and the discussion following it]{goncharovPolylogarithmsRegulatorsArakelov2004}.
\end{remark}

We can consider the following alternative definition of the Grassmannian $n$-logarithm. Let $\mathrm{Herm}^{>0}(\mathbf{C}^{n})$ (respectively $\mathrm{Herm}^{\geq 0}(\mathbf{C}^{n})$) denote the space of positive definite (respectively positive semidefinite) Hermitian forms. Let
\[
  \mathcal{H}_{n} := \mathrm{Herm}^{>0}(\mathbf{C}^{n})/\mathbf{R}^{\times}_{+} \cong \mathrm{SL}_{n}(\mathbf{C})/\mathrm{SU}(n).
\]
The space $\mathcal{H}_{n}$ is a simply connected, non-compact, Riemannian symmetric space of rank $n-1$. There is a naturally defined differential $(2n-1)$-form $\omega_{n}$ on $\mathcal{H}_{n}$ called the \textsf{Borel invariant form}. The form $\omega_{n}$ is obtained from the group trace form $\mathrm{tr}( (g^{-1}\mathrm{d}g)^{2n-1} )$ by restricting at the identity to Hermitian traceless matrices and extending invariantly. If we identify the quotient $\mathcal{H}_{n}$ with determinant-one positive definite Hermitian forms by $H = gg^{*}$, then
\begin{equation}\label{eqn:Borel-form-Hermitian-matrices}
  \omega_{n} = \frac{1}{2^{2n-1}}\mathrm{tr}( (H^{-1}\mathrm{d}H)^{2n-1} ).
\end{equation}
Consider the \textsf{Satake compactification}
\[
  \overline{\mathcal{H}}_{n} := (\mathrm{Herm}^{\geq 0}(\mathbf{C}^{n}) \setminus \{0\})/\mathbf{R}_{+}^{\times}.
\]
The boundary $\partial \overline{\mathcal{H}}_{n}$ is naturally stratified by the rank of the Hermitian form. Let $\mathcal{S}_{r} \subset \overline{\mathcal{H}}_{n}$ denote the stratum of forms of rank $r$. Then $\mathcal{S}_{n} = \mathcal{H}_{n}$ and $\overline{\mathcal{H}}_{n} = \bigsqcup_{r\leq n} \mathcal{S}_{r}$. There is a natural map $\mathcal{S}_{r} \to \mathrm{G}_{n-r}^{r}(\mathbf{C})$ given by $[H] \mapsto \ker H$. Over $\xi \in \mathrm{G}_{n-r}^{r}(\mathbf{C})$ the fibre is $\mathrm{Herm}^{>0}(\mathbf{C}^{n}/\xi)/\mathbf{R}^{\times}_{+}$, so that there is a fibration
\[
  \mathcal{H}_{r} \to \mathcal{S}_{r} \to \mathrm{G}_{n-r}^{r}(\mathbf{C}).
\]
Every hyperplane $h \subset \mathbf{CP}^{n-1}$ defines a positive ray on $\mathrm{Herm}^{\geq 0}(\mathbf{C}^{n})$ given by those positive semidefinite Hermitian forms with kernel $h$. Given a configuration of hyperplanes $h_{0}, \ldots, h_{2n-1} \subset \mathbf{CP}^{n-1}$ with no point in common, denote by $C(h_{0}, \ldots, h_{2n-1})$ the image on $\overline{\mathcal{H}}_{n}$ of the convex hull in $\mathrm{Herm}^{\geq 0}(\mathbf{C}^{n})$ of the rays corresponding to the hyperplanes. We regard this as a possibly degenerate ideal $(2n-1)$-simplex on $\overline{\mathcal{H}}_{n}$, with vertices lying on $\mathcal{S}_{1}$ and with the image of the interior of the parameter simplex lying in $\mathcal{H}_{n}$. Indeed, let $Q_{i} \in \mathrm{Herm}^{\geq 0}(\mathbf{C}^{n})$ with $\mathbf{P}(\ker Q_{i}) = h_{i}$. Then
\begin{equation}\label{eqn:simplex-symmetric-space}
  C(h_{0}, \ldots, h_{2n-1}) = \left\{\left[ \sum_{j=0}^{2n-1} t_{j}Q_{j}\right] \in \overline{\mathcal{H}_{n}} \colon t = (t_{0}, \ldots, t_{2n-1}) \in \Delta^{2n-1}\right\}.
\end{equation}
The simplex $C(h_{0}, \ldots, h_{2n-1})$ does not depend on the choice of $Q_{i}$, although its parametrisation does. Throughout, integration over $C(h_{0}, \ldots, h_{2n-1})$ means integration of the pullback under the parametrisation \eqref{eqn:simplex-symmetric-space}, even when this parametrisation is not injective. One may consider \Cite[\textsection 7.4]{goncharovPolylogarithmsRegulatorsArakelov2004} (where Goncharov attributes the idea to Kontsevich and Nekovar independently) the assignment 
\begin{equation}\label{eqn:nekovar-kontsevich}
  (h_{0}, \ldots, h_{2n-1}) \mapsto \int_{C(h_{0}, \ldots, h_{2n-1})} \omega_{n}.
\end{equation}
The goal of this section is to show that the assignment \eqref{eqn:nekovar-kontsevich} is well-defined (i.e.\ the integral converges) and that it coincides with $\mathcal{L}_{n}^{\mathrm{G}}$ up to a universal constant.
\subsection{Goncharov's cocycle and its infinitesimal form}
Our starting point is that Goncharov \cite{goncharovPolylogarithmsRegulatorsArakelov2004} already made a comparison at the infinitesimal level. Let us unpack what that means. Let $H_{0}, \ldots, H_{2n-1} \in \mathrm{Herm}^{\geq 0}(\mathbf{C}^{n})$ be nonzero. Consider the integral
\begin{equation}\label{eqn:cocycle}
  \psi_{n}(H_{0}, \ldots, H_{2n-1}) := -n^{2n-1}\int_{\mathbf{CP}^{n-1}} \log\frac{H_{1}}{H_{0}} \mathrm{d}\log \frac{H_{2}}{H_{0}} \wedge \cdots \wedge \mathrm{d}\log \frac{H_{2n-1}}{H_{0}},
\end{equation}
where the ratio $\frac{H_{i}}{H_{j}}$ is shorthand for the function on $\mathbf{CP}^{n-1}$ given by
\[
  [v] \mapsto \frac{H_{i}(v,v)}{H_{j}(v,v)}, \quad [v] \in \mathbf{CP}^{n-1}.
\]  
Goncharov \Cite[Prop.\ 5.1, Prop.\ 5.4, Eq.\ (59)]{goncharovPolylogarithmsRegulatorsArakelov2004} shows that: 
\begin{itemize}
  \item \eqref{eqn:cocycle} converges and only depends on the classes of $[H_{i}] \in \overline{\mathcal{H}}_{n}$;
  \item if $H_{0}, \ldots, H_{2n-1}$ are positive definite, then $\psi_{n}$ defines a smooth $\mathrm{SL}_{n}(\mathbf{C})$-invariant $(2n-1)$-cocycle for $\mathcal{H}_{n}$; and
  \item for $2n$ hyperplanes $h_{0}, \ldots, h_{2n-1} \subset \mathbf{CP}^{n-1}$ and rank-one forms $Q_{0}, \ldots, Q_{2n-1}$ with $\mathbf{P}(\ker Q_{i}) = h_{i}$ it follows that
  \begin{equation}\label{eqn:rank-one-cocycle-grassmannian}
    \psi_{n}(Q_{0}, \ldots, Q_{2n-1}) = A_{n}\mathcal{L}^{G}_{n}(h_{0}, \ldots, h_{2n-1}) 
  \end{equation}
  with
  \[
    A_{n} = (-4)^{1-n}(2\pi i)^{n-1}(2n)^{2n-1}\binom{2n-2}{n-1}.
  \]
\end{itemize}
Differentiating the cocycle $\psi_{n}$ along the diagonal via the van Est map
\[
\operatorname{vE}\colon \mathrm{C}^{2n-1}_{\mathrm{sm}}(\mathcal{H}_{n})^{\mathrm{SL}_{n}(\mathbf{C})} \to \Omega^{2n-1}(\mathcal{H}_{n})^{\mathrm{SL}_{n}(\mathbf{C})}
\]
gives a $\mathrm{SL}_{n}(\mathbf{C})$-invariant $(2n-1)$-form $\operatorname{vE}(\psi_{n}) \in \Omega^{2n-1}(\mathcal{H}_{n})^{\mathrm{SL}_{n}(\mathbf{C})}$. Since $\mathcal{H}_{n}$ is a symmetric space, every $\mathrm{SL}_{n}(\mathbf{C})$-invariant form is closed, so $\mathrm{d}(\operatorname{vE}(\psi_{n})) = 0$. Let $[H] \in \mathcal{H}_{n}$ and $K_{1}, \ldots, K_{2n-1} \in \mathrm{T}_{[H]}\mathcal{H}_{n}$. Differentiating the explicit formula \eqref{eqn:cocycle} we obtain (compare with \Cite[Eq. (64)]{goncharovPolylogarithmsRegulatorsArakelov2004} at the identity)
\begin{equation}\label{eqn:infinitesimal-cocycle}
  \mathrm{vE}(\psi_{n})_{[H]}(K_{1}, \ldots, K_{2n-1}) = -n^{2n-1}\int_{\mathbf{CP}^{n-1}} \frac{K_{1}}{H} \mathrm{d}\frac{K_{2}}{H} \wedge \cdots \wedge \mathrm{d}\frac{K_{2n-1}}{H}.
\end{equation}
The key observation is that we can recover $\psi_{n}$ from $\mathrm{vE}(\psi_{n})$ by integrating over the simplex. Indeed, for $H_{0}, \ldots, H_{2n-1} \in \mathrm{Herm}^{>0}(\mathbf{C}^{n})$ consider the simplex
\[
  C(H_{0}, \ldots, H_{2n-1}) := \left\{\left[\sum_{i=0}^{2n-1}t_{i}H_{i}\right] \in \mathcal{H}_{n} \colon t = (t_{0}, \ldots, t_{2n-1}) \in \Delta^{2n-1} \right\} \subset \mathcal{H}_{n}
\]
with vertices at $[H_{0}], \ldots, [H_{2n-1}]$ (observe that a different choice of representatives for $[H_{i}]$ yields the same simplex but with a different parametrisation). 
\begin{theorem}\label{thm:integral-over-algebraic-simplex}
  For $[H_{0}], \ldots, [H_{2n-1}] \in \mathcal{H}_{n}$ we have that
  \[
    \int_{C(H_{0}, \ldots, H_{2n-1})} \mathrm{vE}(\psi_{n}) = \frac{1}{(2n-1)!}\psi_{n}(H_{0}, \ldots, H_{2n-1}).
  \]
\end{theorem}
\begin{proof}
  Follows immediately by applying \cref{thm:barycentric-transgression} to $X = \mathbf{CP}^{n-1}$ and $a_{i}=H_{i}/S$ for some reference positive Hermitian form $S$.
\end{proof}
In \Cite[Theorem 5.7]{goncharovPolylogarithmsRegulatorsArakelov2004} Goncharov shows that (after an appropriate identification in our conventions)
\begin{equation}\label{eqn:van-est-compared-with-Borel}
  \mathrm{vE}(\psi_{n}) = (-1)^{n}\frac{(2n)^{2n-1}(2\pi i)^{n-1}}{(2n-1)(2n-1)!}\omega_{n},
\end{equation}
which will be the tool to relate $\psi_{n}$ with $\omega_{n}$.
\begin{remark}\label{rem:geodesic-vanishing}
  In \Cite[\textsection 7.3]{goncharovPolylogarithmsRegulatorsArakelov2004}, Goncharov asks whether integrating the Borel form over a geodesic simplex recovers $\psi_n$, and states that this is equivalent to its vanishing whenever three vertices lie on the same geodesic. The latter condition fails for $n=3$. Indeed, take
  \[
    H_0=I, \qquad H_1=\operatorname{diag}(2,1,\tfrac12), \qquad H_2=\operatorname{diag}(4,1,\tfrac14).
  \]
  Their classes lie on the geodesic $t\mapsto[\operatorname{diag}(2^t,1,2^{-t})]$, but the matrices are linearly independent. Choose positive definite real symmetric matrices $H_3,H_4,H_5$ completing them to a basis of real symmetric $3\times3$ matrices. After normalising the traces, their convex hull is a nondegenerate five-dimensional simplex in $\mathrm{SL}_3(\mathbf R)/\mathrm{SO}(3)\subset\mathcal H_3$. The restriction of $\omega_3$ to this symmetric space is a nonzero invariant top-degree form, as follows by evaluating the trace formula at the identity, and hence is a volume form. Thus
  \[
    \int_{C(H_0,\ldots,H_5)}\omega_3\neq0.
  \]
  By \cref{thm:integral-over-algebraic-simplex} and \eqref{eqn:van-est-compared-with-Borel}, we obtain $\psi_3(H_0,\ldots,H_5)\neq0$, contradicting the proposed vanishing condition. This gives a negative answer to the geodesic interpretation. The reason why this occurs is because $\mathcal{H}_{3}$ is a symmetric space of rank two, i.e.\ there is a two-dimensional space of invariants of the relative position between two points.
\end{remark}
\subsection{Extension to the rank-one stratum}
In order to prove the comparison we need to extend \cref{thm:integral-over-algebraic-simplex} to the rank-one stratum $\mathcal{S}_{1}$ of $\overline{\mathcal{H}}_{n}$. Throughout this subsection we will keep the following notation. Let $h_{0}, \ldots, h_{2n-1} \subset \mathbf{CP}^{n-1}$ be a configuration of hyperplanes with no point in common and let $Q_{i} \in \mathcal{S}_{1}$ be rank-one positive semidefinite Hermitian forms with $\mathbf{P}(\ker Q_{i}) = h_{i}$. Let $\Phi\colon \Delta^{2n-1} \to \overline{\mathcal{H}}_{n}$ be the parametrisation 
\[
  t \mapsto \left[ \sum_{i=0}^{2n-1} t_{i}Q_{i}\right]
\]
of $C(Q_{0}, \ldots, Q_{2n-1})$.
Throughout this subsection, integration over $C(Q_{0}, \ldots, Q_{2n-1})$ means integration of the pullback under $\Phi$, including when the parametrisation is not injective. Write (which is the notation used locally in \cref{thm:barycentric-transgression})
\[
  S:= \sum_{i=0}^{2n-1}Q_{i}, \quad a_{i}:=\frac{Q_{i}}{S}, \quad b_{i}:=a_{i}-a_{0}, \quad A_{t}:=\sum_{i=0}^{2n-1}t_{i}a_{i},
\]
for $t = (t_{0}, \ldots, t_{2n-1}) \in \Delta^{2n-1}$. Observe that $S$ is positive definite. Indeed, since the hyperplanes $h_{0}, \ldots, h_{2n-1}$ have no point in common, we have that
\[
  \ker S = \bigcap_{i=0}^{2n-1}\ker Q_{i} = 0.
\]
In particular, the functions $a_{i}$ are smooth nonnegative functions on $\mathbf{CP}^{n-1}$. Moreover, for every $t \in \mathring{\Delta}^{2n-1}$, $A_{t}$ is strictly positive.
\begin{proposition}\label{prop:absolute-integrability}
    The pullback $\Phi^{*}(\mathrm{vE}(\psi_{n}))$ is absolutely integrable on $\mathring{\Delta}^{2n-1}$. In particular, from \eqref{eqn:van-est-compared-with-Borel}, so is $\omega_{n}$.
\end{proposition}
\begin{proof}
  The idea is to use an adaptation of the proof of \cref{thm:barycentric-transgression} to obtain an estimate of the integral. From the explicit formula \eqref{eqn:infinitesimal-cocycle} of $\mathrm{vE}(\psi_{n})$, its pullback $\Phi^{*}\mathrm{vE}(\psi_{n})$ is given by
  \[
    -n^{2n-1}\left(\int_{\mathbf{CP}^{n-1}} \frac{b_{1}}{A_{t}} \mathrm{d}\frac{b_{2}}{A_{t}} \wedge \cdots \wedge \mathrm{d}\frac{b_{2n-1}}{A_{t}}\right)\,\mathrm{d}t_{1}\wedge \cdots \mathrm{d}t_{2n-1}.
  \]
  Recall from the proof of \cref{thm:barycentric-transgression} the radial form $\rho(b)$ defined by the formula

\[
    \rho(b)=\sum_{i=1}^{2n-1}(-1)^{i+1}b_{i}\mathrm{d}b_{1}\wedge \cdots \wedge \widehat{\mathrm{d}b_{i}} \wedge \cdots \wedge \mathrm{d}b_{2n-1}.
  \]
For each $t \in \mathring{\Delta}^{2n-1}$, use Stokes' theorem and the formula
  \[
    \rho\left( \frac{b}{A_{t}} \right) = \frac{\rho(b)}{A_{t}^{2n-1}}
  \]
  to rewrite $\Phi^{*}\mathrm{vE}(\psi_{n})$ as
  \[
    - \frac{n^{2n-1}}{2n-1}\left(\int_{\mathbf{CP}^{n-1}} \frac{\rho(b)}{A_{t}^{2n-1}}\right)\,\mathrm{d}t_{1}\wedge \cdots \mathrm{d}t_{2n-1}.
  \]
  Therefore, it suffices, from the triangle inequality and Tonelli's theorem, to show that
  \[
    \int_{\mathbf{CP}^{n-1}}|\rho(b)|\left(\int_{\Delta^{2n-1}}A_{t}^{-(2n-1)}\mathrm{d}t\right) < \infty.
  \]
  As in the proof of \cref{thm:barycentric-transgression}, by the Genocchi-Hermite interpolation formula for the logarithm \eqref{eqn:Hermite-Genocchi-logarithm} we have that
  \[
    \int_{\Delta^{2n-1}} A_{t}^{-(2n-1)}\mathrm{d}t = \frac{1}{(2n-2)!}[a_{0}, \ldots, a_{2n-1}]\log.
  \]
  Now, we use the integral representation of the logarithm
  \[
    \log x = \int_{0}^{\infty} \left( \frac{1}{1+s} - \frac{1}{x+s}\right)\mathrm{d}s.
  \]
  Observe that
  \[
    [a_{0}, \ldots, a_{2n-1}] \frac{1}{x+s} = - \frac{1}{\prod_{i=0}^{2n-1}(a_{i}+s)}.
  \]
  Therefore,
  \[
    \int_{\Delta^{2n-1}}A_{t}^{-(2n-1)} \mathrm{d}t = \frac{1}{(2n-2)!}\int_{0}^{\infty}\frac{\mathrm{d}s}{\prod_{i=0}^{2n-1}(a_{i}+s)}.
  \]
  Expanding $\rho(b)$ using $b_{i}=a_{i}-a_{0}$, we obtain 

  \[
    \rho(b)=\sum_{0\leq p<q\leq 2n-1}(-1)^{p+q}(a_p-a_q)
    \bigwedge_{i\neq p,q}\mathrm{d}a_i,
  \]
  Away from the hyperplanes, we have
\[
    \begin{aligned}
    &|a_p-a_q|\left(\prod_{i\neq p,q}a_i\right)
      \int_0^\infty\frac{\mathrm{d}s}{\prod_i(a_i+s)}\\
    &\qquad\leq |a_p-a_q|\int_0^\infty
      \frac{\mathrm{d}s}{(a_p+s)(a_q+s)}
      =\left|\log\frac{a_p}{a_q}\right|.
    \end{aligned}
  \]
  The inequality follows from $a_i/(a_i+s)\leq 1$. The last equality follows by partial fractions when $a_p\neq a_q$. When $a_p=a_q$, both sides vanish. Using $\mathrm{d}a_i=a_i\,\mathrm{d}\log a_i$ in the expansion of $\rho(b)$ and applying the triangle inequality gives
\[
    |\rho(b)|\left(\int_{\Delta^{2n-1}}A_t^{-(2n-1)}\mathrm{d}t\right)
    \leq\frac{1}{(2n-2)!}\sum_{p<q}
      \left|\log\frac{a_p}{a_q}\right|
      \left|\bigwedge_{i\neq p,q}\mathrm{d}\log a_i\right|
  \]
  relative to a smooth volume density. The remaining part of the argument is essentially that of \Cite[Lemma 2.5]{goncharovPolylogarithmsRegulatorsArakelov2004}. Resolve the union of the hyperplanes by a projective resolution with total transform with simple normal crossings. The preceding inequality also holds after pullback. Fix a point $P$ on the resolution and choose holomorphic coordinates $z_j$ in which the components of the total transform passing through $P$ are coordinate hyperplanes. We continue to write $a_i$ for their pullbacks. Locally, each $a_i$ is a positive smooth factor times a monomial in the $|z_j|^2$. Therefore, each $\log a_i$ is a smooth function plus a linear combination of the $\log|z_j|$, and each $\mathrm{d}\log a_i$ is a smooth one-form plus a linear combination of the $\mathrm{d}\log|z_j|$. Expanding the wedges we see that the terms containing two copies of the same $\mathrm{d}\log|z_j|$ vanish. Shrinking the polydisk so that $|z_j|<1$, we see that, relative to the Euclidean volume density on the resolution, we can bound by
\[
      \mathrm{constant}\times
      \left(1+\sum_j|\log|z_j||\right)\prod_j|z_j|^{-1},
   \]
   where $j$ ranges through the hyperplanes $\{z_j=0\}$ of the total transform passing through $P$. Since
   \[
  \int_{|z_{j}| < \varepsilon} \frac{\mathrm{d}(\mathfrak{Re}z_{j}) \mathrm{d}(\mathfrak{Im}z_{j})}{|z_{j}|} = \int_{0}^{2\pi}\int_{0}^{\varepsilon} \mathrm{d}r_{j}\mathrm{d}\theta_{j} = 2\pi \varepsilon,
   \]
   all the factors $|z_{j}|^{-1}$ are integrable. The logarithmic factors are also integrable, since $\int_0^1|\log r_j|\,\mathrm{d}r_j<\infty$. Fubini gives that the products of the factors are integrable on a polydisk. Also, any smooth volume density is bounded by a constant times the Euclidean volume on the polydisk. Since $P$ was arbitrary and the resolution is compact, doing this in a finite collection of polydisks gives global integrability, as desired. The resolution is an isomorphism off the hyperplanes, so this also proves integrability in $\mathbf{CP}^{n-1}$.
\end{proof}
We now approximate the rank-one forms $Q_i$ by positive definite forms and show that the corresponding values of $\psi_n$ converge.
\begin{proposition}\label{prop:convergence-under-regularisation}
  Set $Q_{i}^{\varepsilon} = Q_{i} + \varepsilon S \in \mathrm{Herm}^{> 0 }(\mathbf{C}^{n})$. Then
  \[
    \lim_{\varepsilon \to 0^{+}} \psi_{n}(Q_{0}^{\varepsilon}, \ldots, Q_{2n-1}^{\varepsilon}) = \psi_{n}(Q_{0}, \ldots, Q_{2n-1}).
  \]
\end{proposition}
\begin{proof}
  We make the local argument of \Cite[Lemma 2.5]{goncharovPolylogarithmsRegulatorsArakelov2004} uniform in $\varepsilon$, so that we can compute the limit via dominated convergence. We proceed in the same way as in the last part of the proof of \cref{prop:absolute-integrability}. Fix a point $P$ on the same resolution and choose coordinates $z_j$ in a polydisk around $P$ so that the components of the total transform passing through $P$ are distinct coordinate hyperplanes. Again, write $a_i$ for their pullbacks. Shrink the polydisk so that $|z_j|<1$. Each $a_i$ is a positive smooth factor times a monomial in the $|z_j|^2$. First, we obtain bounds independent of $0 < \varepsilon \leq 1$. In this case, $0 \leq a_{i} \leq 1$, so
  \[
    a_{i} \leq a_{i} + \varepsilon \leq 2.
  \]
  Therefore, away from the total transform, $|\log(a_i+\varepsilon)|\leq |\log a_i|+\log 2$, and the local monomial expression gives
  \[
    |\log(a_i+\varepsilon)|\leq\mathrm{constant}\times
    \left(1+\sum_j|\log|z_j||\right).
  \]
In addition,
  \[
    \mathrm{d}\log(a_i+\varepsilon)=\frac{a_i}{a_i+\varepsilon}\,\mathrm{d}\log a_i,
    \qquad 0\leq\frac{a_i}{a_i+\varepsilon}\leq 1.
  \]
The scalar factor is bounded independently of $\varepsilon$. Each $\mathrm{d}\log a_i$ is a smooth one-form plus a linear combination of the $\mathrm{d}\log|z_j|$. Now, using the formula \eqref{eqn:cocycle}, we can write down the pullback of the integrand defining $\psi_{n}(Q_{0}^{\varepsilon}, \ldots, Q_{2n-1}^{\varepsilon})$ as
  \[
    -n^{2n-1}\log \frac{a_{1} + \varepsilon}{a_{0} + \varepsilon} \mathrm{d}\log \frac{a_{2} + \varepsilon}{a_{0} + \varepsilon} \wedge \cdots \wedge \mathrm{d}\log \frac{a_{2n-1} + \varepsilon}{a_{0} + \varepsilon}.
  \]
  If we expand out the wedge, each differential form with index $i \geq 2$ appears at most once. The terms that contain two copies of $\mathrm{d}\log(a_{0} + \varepsilon)$ vanish. The terms containing two copies of the same $\mathrm{d}\log|z_j|$ also vanish. Therefore, each coordinate contributes at most one factor $|z_j|^{-1}$. Together with the logarithm bound, we have that the integrand of $\psi_{n}(Q_{0}^{\varepsilon}, \ldots, Q_{2n-1}^{\varepsilon})$ is bounded by
  \[
    \mathrm{constant} \times \left( 1 + \sum_{j}|\log|z_{j}|| \right) \prod_{j} |z_{j}|^{-1},
  \]
  where the constant is independent of $\varepsilon$, relative to the Euclidean volume on the resolution. This bound is integrable by the same polar-coordinate calculation used in the proof of \cref{prop:absolute-integrability}. The coordinates that do not define hyperplanes contribute only finite volume factors. Therefore, since $P$ was arbitrary, by compactness of the resolution, we can do this in only finitely many polydisks. The pulled-back integrands converge pointwise off the total transform, and integrals in the resolution agree with those in $\mathbf{CP}^{n-1}$ off the divisors. Applying the dominated convergence theorem (after some partition of unity in each polydisk) establishes the desired result.
\end{proof}
We can now extend \cref{thm:integral-over-algebraic-simplex} to the stratum $\mathcal{S}_{1}$.
\begin{theorem}\label{thm:integration-extends-to-rank-one-stratum}
  For every configuration of hyperplanes $h_{0}, \ldots, h_{2n-1} \subset \mathbf{CP}^{n-1}$ with no point in common, and rank one Hermitian forms $Q_{i} \in \mathcal{S}_{1}$ with $\mathbf{P}(\ker Q_{i}) = h_{i}$, we have that
  \[
    \int_{C(h_{0}, \ldots, h_{2n-1})} \mathrm{vE}(\psi_{n}) = \frac{1}{(2n-1)!}\psi_{n}(Q_{0}, \ldots, Q_{2n-1}).
  \]
\end{theorem}
\begin{proof}
  Consider, as in \cref{prop:convergence-under-regularisation}, the regularised forms $Q_{i}^{\varepsilon} = Q_{i} + \varepsilon S$ which are positive. From \cref{thm:integral-over-algebraic-simplex} we have that
  \begin{equation}\label{eqn:regularised-inetgral}
    \int_{C(Q_{0}^{\varepsilon}, \ldots, Q_{2n-1}^{\varepsilon})} \mathrm{vE}(\psi_{n}) = \frac{1}{(2n-1)!}\psi_{n}(Q_{0}^{\varepsilon}, \ldots, Q_{2n-1}^{\varepsilon}).
  \end{equation}
  Observe that
  \[
    \left[ \sum_{i=0}^{2n-1} t_{i}Q_{i}^{\varepsilon} \right] = \left[ \sum_{i=0}^{2n-1} \frac{t_{i} + \varepsilon}{1 + 2n \varepsilon} Q_{i} \right].
  \]
  The orientation-preserving transformation
  \[
    t \mapsto \left( \frac{t_{0} + \varepsilon}{1 + 2n \varepsilon}, \ldots, \frac{t_{2n-1} + \varepsilon}{1 + 2n \varepsilon} \right)
  \]
  identifies $\Delta^{2n-1}$ with
  \[
  \left\{ t \in \Delta^{2n-1} \colon t_{i} \geq \frac{\varepsilon}{1+2n \varepsilon} \text{ for all }i  \right\}.
  \]
  Hence, the left-hand side of \eqref{eqn:regularised-inetgral} is the integral over the original simplex pulled back to the truncated domain. In the limit $\varepsilon \to 0^{+}$, these truncated domains exhaust the simplex. Therefore, from absolute integrability obtained in \cref{prop:absolute-integrability} gives that
  \[
    \lim_{\varepsilon \to 0^{+}} \int_{C(Q_{0}^{\varepsilon}, \ldots, Q_{2n-1}^{\varepsilon})} \mathrm{vE}(\psi_{n}) = \int_{C(h_{0}, \ldots, h_{2n-1})} \mathrm{vE}(\psi_{n}).
  \]
  From \cref{prop:convergence-under-regularisation} the right-hand side of \eqref{eqn:regularised-inetgral} converges to $1/(2n-1)!\psi_{n}(Q_{0}, \ldots, Q_{2n-1})$, so we are done. 
\end{proof}
\begin{corollary}\label{cor:integral-Borel-equals-grassmannian}
  For a configuration of hyperplanes $h_{0}, \ldots, h_{2n-1} \subset \mathbf{CP}^{n-1}$ with no point in common we have that
  \[
    \int_{C(h_{0}, \ldots, h_{2n-1})} \omega_{n} = C_{n} \mathcal{L}_{n}^{\mathrm{G}}(h_{0}, \ldots, h_{2n-1}),
  \]
  where
  \[
    C_{n} = - \frac{2n-1}{4^{n-1}}\binom{2n-2}{n-1}.
  \]
\end{corollary}
\begin{proof}
  Follows from \eqref{eqn:rank-one-cocycle-grassmannian}, \eqref{eqn:van-est-compared-with-Borel}, and \cref{thm:integration-extends-to-rank-one-stratum}.
\end{proof}
\subsection{A formula of Dan}
Let $v_{0}, \ldots, v_{2n-1} \in \mathbf{C}^{n}$ be a generic configuration of nonzero vectors. To each $v_{i}$ we associate the rank-one Hermitian matrix $H_{i} = v_{i}v_{i}^{*}$ and the hyperplane $h_{i} = \mathbf{P}(\ker(v_{i}^{*}))$. We write $\mathcal{L}_{n}^{\mathrm{G}}(v_{0}, \ldots, v_{2n-1})$ for $\mathcal{L}_{n}^{\mathrm{G}}(h_{0}, \ldots, h_{2n-1})$. The corresponding Hermitian simplex is parametrised by $H_{t} = \sum_{i=0}^{2n-1}t_{i}v_{i}v_{i}^{*}$. We will expand the integral in \cref{cor:integral-Borel-equals-grassmannian} explicitly. Since $\mathrm{d}H_{t} = \sum_{i=0}^{2n-1}v_{i}v_{i}^{*}\mathrm{d}t_{i}$, multiplication of rank-one matrices gives
\[
  \mathrm{tr}\bigl((H_{t}^{-1}\mathrm{d}H_{t})^{2n-1}\bigr)
  = \sum_{i_{1},\ldots,i_{2n-1}=0}^{2n-1}
  \left(\prod_{a=1}^{2n-1}v_{i_{a}}^{*}H_{t}^{-1}v_{i_{a+1}}\right)
  \mathrm{d}t_{i_{1}}\wedge\cdots\wedge\mathrm{d}t_{i_{2n-1}},
\]
where $i_{2n}=i_{1}$. The trace expression is unchanged by multiplication of $H_t$ by a positive scalar function, so \eqref{eqn:Borel-form-Hermitian-matrices} applies without imposing $\det(H_t)=1$. Hence
\[
  \int_{\Delta^{2n-1}}\mathrm{tr}\bigl((H_t^{-1}\mathrm{d}H_t)^{2n-1}\bigr)
  = 2^{2n-1}\int_{C(h_0,\ldots,h_{2n-1})}\omega_n.
\]
For $I \subset \{0, \ldots, 2n-1\}$ write
\[
  t_{I} = \prod_{i\in I}t_{i}, \quad v_{I}=(v_{i})_{i\in I},
  \quad \Delta_{I}=\mathrm{det}(v_{I})\quad\text{if }|I|=n,
\]
with vectors ordered increasingly. By Cauchy-Binet we get
\[
  \det(H_{t}) = \sum_{|I|=n}t_{I}|\Delta_{I}|^{2}
\]
and similarly \[v_{i}^{*}\mathrm{adj}(H_{t})v_{j} = \sum_{|I|=n-1}t_{I}\overline{\det(v_{I}, v_{i})}\det(v_{I}, v_{j}).\]
Therefore,
\[
  v_{i}^{*}H_{t}^{-1}v_{j} = \frac{\sum_{|I|=n-1}t_{I}\overline{\det(v_{I}, v_{i})}\det(v_{I}, v_{j})}{\sum_{|I|=n}t_{I}|\Delta_{I}|^{2}}.
\]
Substituting these expressions into the trace expansion and applying \cref{cor:integral-Borel-equals-grassmannian} gives the following determinant formula.
\begin{theorem}\label{thm:president-formula}
For a generic configuration of vectors $v_{0}, \ldots, v_{2n-1}$ in $\mathbf{C}^{n}$ we have that $c_{n}\mathcal{L}^{\mathrm{G}}_{n}(v_{0}, \ldots, v_{2n-1})$ equals
\[
  \int_{\Delta^{2n-1}}
  \frac{\displaystyle
    \sum_{\substack{|I_{1}|=\cdots=|I_{2n-1}|=n-1 \\ 0\leq i_{1},\ldots,i_{2n-1}\leq 2n-1}}
    \left(\prod_{a=1}^{2n-1}t_{I_a}\overline{\det(v_{I_a},v_{i_a})}\det(v_{I_a},v_{i_{a+1}})\right)
    \mathrm{d}t_{i_1}\wedge\cdots\wedge\mathrm{d}t_{i_{2n-1}}
  }{\displaystyle\left(\sum_{|I|=n}t_I|\Delta_I|^2\right)^{2n-1}},
\]
where $i_{2n}=i_{1}$ and
\[
  c_{n}=-2(2n-1)\binom{2n-2}{n-1}.
\]
\end{theorem}
The determinant integral is the one appearing in \Cite[Thm. 4.1, Eq. (6)]{danMultitransgressionRegulators2008}. Including Dan's prefactor $-(n-1)!/(2(2n-1)!)$, the right-hand side of his formula equals
\[
  \frac{1}{(n-1)!}\mathcal{L}^{\mathrm{G}}_{n}(h_0,\ldots,h_{2n-1}).
\]
This identifies Dan's explicit integral with the Grassmannian polylogarithm which establishes part of his conjecture \Cite[Conj. 4.2]{danMultitransgressionRegulators2008}. The further identification with the maximal period of the mixed Hodge structure is not established here.
\subsection{Canonical integrals}\label{sec:canonical-integrals}
The content of this subsection is an application of \cref{cor:integral-Borel-equals-grassmannian} in the context of \textsf{canonical integrals}, which are generalised Feynman integrals of the kind studied in quantum field theory. We will recover the convergence of primitive canonical integrals that appeared in \cite{brownInvariantDifferentialForms2021} and the evaluation of wheel integrals of \cite{brownWheelClassesLocally2025}.

We sketch the setting of \cite{brownWheelClassesLocally2025}. Let $G$ be a finite connected graph with edge set $E(G)$. Write $h(G):= \mathrm{rank}(\mathrm{H}_{1}(G;\mathbf{Z}))$. Choose an orientation on each edge, which identifies the group of $1$-chains with $\mathbf{Z}^{E(G)}$ and gives $\mathrm{H}_{1}(G;\mathbf{Z}) \subset \mathbf{Z}^{E(G)}$ as the subgroup of cycles. Associate a variable $x_{e}$ to each edge $e \in E(G)$. Endow $E(G)$ with the pairing $\langle e,e'\rangle = \delta_{e,e'}x_{e}$. Choose an integral cycle basis $c_{1}, \ldots, c_{h(G)}$ of $\mathrm{H}_{1}(G;\mathbf{Z})$. The \textsf{graph Laplacian} is given by
\begin{equation}\label{eqn:graph-laplacian}
  \Lambda_{G} := (\langle c_{i}, c_{j}\rangle)_{1\leq i,j \leq h(G)}.
\end{equation}
It is a symmetric, positive-definite matrix for all $x_{e} > 0$ \Cite[\textsection 1.1]{brownWheelClassesLocally2025}. It is also independent of the choice of edge directions. Changing the cycle basis replaces $\Lambda_{G}$ by $P^{T}\Lambda_{G}P$ for some $P \in \mathrm{GL}_{h(G)}(\mathbf{Z})$. For a symmetric matrix $X$, Brown-Schnetz \cite{brownWheelClassesLocally2025} consider
\[
  \omega_{X}^{k} := \mathrm{tr}( (X^{-1}\mathrm{d}X)^{k} ).
\]
Now, let $n > 1$ odd and assume that the graph has $2n$ edges. The \textsf{canonical integral} is given by
\[
  I_{G}(\omega^{2n-1}) := \int_{\sigma_{G}} \omega^{2n-1}_{\Lambda_{G}},
\]
where
\[
  \sigma_{G} := \{[x_{e}]_{e}\colon x_{e} \geq 0\} \subset \mathbf{RP}^{|E(G)|-1}
\]
i.e.\ $\Delta^{2n-1}$. Brown shows \Cite[Corollary 8.3 and Lemma 8.4]{brownInvariantDifferentialForms2021} that the canonical integral $I_{G}(\omega^{2n-1})$ converges; it vanishes unless $G$ is bridgeless and $|E(G)|=2h(G)$. Henceforth, let us assume that $G$ is bridgeless, $|E(G)|=2n$, and $h(G) = n$. The chosen cycle basis $c_{1}, \ldots, c_{n}$ identifies $\mathrm{H}_{1}(G; \mathbf{C})$ with $\mathbf{C}^{n}$. For each edge $e$, consider the Hermitian form
\[
  Q_{e}(z_{1}, \ldots, z_{n}) := \left| \sum_{i=1}^{n} (c_{i})_{e} z_{i} \right|^{2},
\]
where $(c_{i})_{e}$ is the coefficient of $e$ in the cycle $c_{i}$, as in the definition of $\Lambda_{G}$. The associated hyperplane is
\[
  h_{e}:= \left\{[z_{1}: \cdots : z_{n}] \in \mathbf{CP}^{n-1} \colon \sum_{i=1}^{n} (c_{i})_{e} z_{i} = 0\right\}.
\]
Since $G$ is bridgeless, every edge belongs to a cycle. Thus, for each $e$, the coefficients $(c_{i})_{e}$ are not all zero and $Q_{e}$ is rank-one. Moreover, $\bigcap_{e \in E(G)} h_{e} = \emptyset$ because $c_{i}$ forms a basis. From the definition of graph Laplacian, we see that
\[
  \Lambda_{G}(x) = \sum_{e \in E(G)}x_{e}Q_{e}.
\]
Fix an ordering $E(G) = \{e_{1}, \ldots, e_{2n}\}$.
In the comparison below, $\sigma_{G}$ is endowed with the simplex orientation determined by this edge ordering.
This ordering identifies
\[
  [x_{e}] \mapsto [\Lambda_{G}(x)]
\]
with the parametrisation of $C(h_{e_{1}}, \ldots, h_{e_{2n}})$. We are ready to prove a comparison of the canonical integral $I_{G}(\omega^{2n-1})$ with the Grassmannian $n$-logarithm $\mathcal{L}^{\mathrm{G}}_{n}(h_{e_{1}}, \ldots, h_{e_{2n}})$.
\begin{theorem}\label{thm:borwn-schnetz-comparison}
  Let $G$ be a connected bridgeless graph with $2n$ ordered edges $e_{1}, \ldots, e_{2n}$ and $h(G) = n$. Then
  \[
    I_{G}(\omega^{2n-1}) = -n\binom{2n}{n}\mathcal{L}_{n}^{\mathrm{G}}(h_{e_{1}}, \ldots, h_{e_{2n}}).
  \]
\end{theorem}
\begin{proof}
  By the previous discussion,
  \[
    \Phi\colon \sigma_{G} \to \overline{\mathcal{H}}_{n}, \quad [x_{e}] \mapsto [\Lambda_{G}(x)] = \left[ \sum_{e}x_{e}Q_{e} \right]
  \]
  parametrises $C(h_{e_{1}}, \ldots, h_{e_{2n}})$. The hyperplanes have no point in common, so \cref{prop:absolute-integrability} and \cref{cor:integral-Borel-equals-grassmannian} apply. In our normalisation
  \[
    \Phi^{*}\omega_{n} = 2^{-(2n-1)}\omega_{\Lambda_{G}}^{2n-1},
  \]
  so that
  \[
    I_{G}(\omega^{2n-1}) = 2^{2n-1}\int_{C(h_{e_{1}}, \ldots, h_{e_{2n}})} \omega_{n} = 2^{2n-1}C_{n}\mathcal{L}_{n}^{\mathrm{G}}(h_{e_{1}}, \ldots, h_{e_{2n}}) = -n\binom{2n}{n}\mathcal{L}_{n}^{\mathrm{G}}(h_{e_{1}}, \ldots, h_{e_{2n}}).
  \]
  Also, in this case \cref{prop:absolute-integrability} gives the convergence of the left-hand side which was previously established in \Cite[Corollary 8.3]{brownInvariantDifferentialForms2021}.
\end{proof}
As a corollary, we recover the properties of canonical integrals under edge contraction and edge deletion from the known properties of the Grassmannian $n$-logarithm in \cref{thm:properties-grassmannian-polylogarithm}.
\begin{corollary}[{\Cite[\textsection 8.2]{brownInvariantDifferentialForms2021}}]\label{cor:contraction-deletion}
\leavevmode
\begin{enumerate}
  \item Let $G$ be a connected bridgeless graph without tadpoles, with $|E(G)| = 2n+1$ and $h(G)=n$. Then
  \[
    \sum_{i=1}^{2n+1}(-1)^{i}I_{G/e_{i}}(\omega^{2n-1}) = 0.
  \]
  \item Let $G$ be a graph without tadpoles with $|E(G)|=2n+1$ and $h(G)=n+1$, and suppose that deleting any edge leaves a connected bridgeless graph. Then
  \[
    \sum_{i=1}^{2n+1}(-1)^{i}I_{G\backslash e_{i}}(\omega^{2n-1}) = 0.
  \]
\end{enumerate}
\end{corollary}
To recover the main theorem of \cite{brownWheelClassesLocally2025} we now focus on wheel graphs. For $n \geq 3$, the \textsf{wheel graph} $W_{n}$ is obtained from a cycle with $n$ vertices by adjoining one central vertex and joining it to every vertex of the cycle (see \Cite[Figure 1]{brownWheelClassesLocally2025}). The original cycle is called the \textsf{rim}, and the $n$ edges joining it to the central vertex are called \textsf{spokes}. Since $W_{n}$ has $n+1$ vertices and $2n$ edges and it is connected, we must have $h(W_{n}) = n$. Also, the wheel graph is bridgeless, hence it satisfies the hypotheses of \cref{thm:borwn-schnetz-comparison}. In the usual planar drawing of $W_{n}$, the wheel has $n$ triangular faces. The spokes form a spanning tree. Any rim edge determines precisely a triangle containing that edge and its two adjacent spokes. Therefore, the oriented boundaries of these triangles form an integral cycle basis. Number the rim edges $e_{1}, \ldots, e_{n}$ anticlockwise and orient them in that direction, with the spoke at its endpoint inward and the spoke at its starting point outward. The last spoke is traversed against its orientation, so its boundary cycle is
\[
  c_{1} = e_{1} + e_{n+1} - e_{2n}.
\]
For the rest, we have
\[
  c_{i} = e_{i} + e_{n+i} - e_{n+i-1}, \quad 2 \leq i \leq n.
\]
They form an integral basis of $\mathrm{H}_{1}(W_{n};\mathbf{Z})$. In this basis
\[
  Q_{e_{i}}(z_{1}, \ldots, z_{n}) = |z_{i}|^{2}, \quad 1 \leq i \leq n.
\]
For $1 \leq i < n$, the spoke $e_{n+i}$ occurs with coefficient $+1$ in $c_{i}$ and $-1$ in $c_{i+1}$. The last spoke occurs with coefficient $+1$ in $c_{n}$ and $-1$ in $c_{1}$. Hence,
\[
  Q_{e_{n+i}}(z) = |z_{i} - z_{i+1}|^{2} \quad (1 \leq i < n), \quad Q_{e_{2n}}(z) = |z_{n}-z_{1}|^{2}.
\]
Therefore, the ordered hyperplane configuration associated with $W_{n}$ is
\begin{equation}\label{eqn:special-hyperplane-configuration}
  \{z_{1}=0\}, \ldots, \{z_{n}=0\}, \{z_{1}-z_{2}=0\}, \ldots, \{z_{n-1} - z_{n} = 0\}, \{z_{n}-z_{1}=0\}.
\end{equation}
The configuration of hyperplanes \eqref{eqn:special-hyperplane-configuration} is special in the sense of Goncharov. A \textsf{special configuration} \Cite[\textsection 4.4]{goncharovPolylogarithmsRegulatorsArakelov2004} is a configuration of $2n$ points $(l_{1}, \ldots, l_{n}, m_{1}, \ldots, m_{n})$ in $\mathbf{CP}^{n-1}$, where $l_{1}, \ldots, l_{n}$ are vertices of a simplex and $m_{i}$ lies on the line joining $l_{i}$ to $l_{i+1}$ away from its endpoints. Indices are cyclic. The moduli of such special configurations is $\mathbf{P}^{1}\setminus\{0, \infty\}$ via its \textsf{generalised cross-ratio} (see \Cite[Proposition 4.4]{goncharovPolylogarithmsRegulatorsArakelov2004} for details). Goncharov explicitly computes the values of the Grassmannian $n$-logarithm at special configurations.
\begin{theorem}[{\Cite[Theorem 4.5]{goncharovPolylogarithmsRegulatorsArakelov2004}}]\label{thm:special-configurations-evaluation}
  For a special configuration $(l_{1}, \ldots, l_{n}, m_{1}, \ldots, m_{n})$ with generalised cross-ratio $s$ the value of $\mathcal{L}_{n}^{\mathrm{G}}$ is
  \[
    -(-1)^{n(n-1)/2}4^{n-1}\binom{2n-2}{n-1}^{-1}\tilde{\mathcal{L}}_{n}(s),
  \]
  where $\tilde{\mathcal{L}}_{n}$ denotes Levin's modified single-value polylogarithm.
\end{theorem}
We only need the value at $s=1$. For odd $n$,\footnote{The constant in \cite{goncharovPolylogarithmsRegulatorsArakelov2004} contains a typo, we have corrected it with the version from \Cite[\textsection 4.2]{levinNoteRHodgeTateSheaves2001}.}
\[
  \tilde{\mathcal{L}}_{n}(1) = \frac{(2n-3)!!}{(2n-2)!!}\zeta(n) = 4^{1-n}\binom{2n-2}{n-1}\zeta(n).
\]
Therefore, the value of $\mathcal{L}_{n}^{\mathrm{G}}$ on a special configuration of generalised cross-ratio $1$ is
\[
  (-1)^{n(n-1)/2 + 1}\zeta(n).
\]
\begin{proposition}\label{prop:value-grassmannian-special-configuration}
  For the configuration \eqref{eqn:special-hyperplane-configuration} and odd $n \geq 3$,
  \[
    \mathcal{L}_{n}^{\mathrm{G}}(h_{e_{1}}, \ldots, h_{e_{2n}}) = (-1)^{n(n-1)/2}\zeta(n).
  \]
\end{proposition}
\begin{proof}
  Consider the ordered points
  \[
    [z_{1}], \ldots, [z_{n}], [z_{1}-z_{2}], \ldots, [z_{n-1}-z_{n}],[z_{n}-z_{1}]
  \]
  in the projective space of linear forms on $\mathbf{C}^{n}$. They clearly form a special configuration. In the quotient
  \[
    (\mathbf{C}z_{i} \oplus \mathbf{C}z_{i+1})/\mathbf{C}(z_{i}-z_{i+1})
  \]
  the images of $z_{i}$ and $z_{i+1}$ coincide. The successive identifications send
  \[
    z_{1} \mapsto z_{2} \mapsto \cdots \mapsto z_{n} \mapsto z_{1}
  \]
  and the composition is the identity. Hence, their generalised cross-ratio is $1$ (see \cite[Proposition 4.4]{goncharovPolylogarithmsRegulatorsArakelov2004}). 
  To compare the orderings, let $D$ be the matrix whose $i$th column gives the coefficients of $z_{i}-z_{i+1}$ in the basis $z_{1},\ldots,z_{n}$, with indices taken modulo $n$. The special configuration is represented by $(I\mid D)$. The correspondence in \Cite[\textsection 4.1]{goncharovPolylogarithmsRegulatorsArakelov2004} uses the kernel of this matrix. A basis of the kernel gives the matrix of hyperplane equations $(-D^{\mathsf{t}}\mid I)$. Swapping the two blocks gives $(I\mid -D^{\mathsf{t}})$, with sign $(-1)^{n^{2}}$. The columns of $-D^{\mathsf{t}}$ are the columns of $D$ in the order $n,1,\ldots,n-1$. A cyclic permutation of these $n$ columns recovers the wheel ordering, with sign $(-1)^{n-1}$. The total sign is therefore
  \[
    (-1)^{n^{2}}(-1)^{n-1}=-1.
  \]
  Therefore,
  \[
    \mathcal{L}_{n}^{\mathrm{G}}(h_{e_{1}}, \ldots, h_{e_{2n}}) = -(-1)^{n(n-1)/2 + 1}\zeta(n) = (-1)^{n(n-1)/2}\zeta(n).
  \]
\end{proof}
\begin{theorem}[{\Cite[Theorem 1]{brownWheelClassesLocally2025}}]\label{thm:evaluation-wheel-classes}
  For all odd $n \geq 3$, with Brown-Schnetz's orientation,
  \[
    I_{W_{n}}(\omega^{2n-1}) = n \binom{2n}{n}\zeta(n).
  \]
\end{theorem}
\begin{proof}
  Using the ordered simplex orientation corresponding to $(h_{e_{1}}, \ldots, h_{e_{2n}})$, from \cref{thm:borwn-schnetz-comparison} and \cref{prop:value-grassmannian-special-configuration} we obtain
    \begin{align*}
      \int_{\sigma_{W_{n}}} \omega_{\Lambda_{W_{n}}}^{2n-1} &= -n  \binom{2n}{n} \mathcal{L}_{n}^{\mathrm{G}}(h_{e_{1}}, \ldots h_{e_{2n}}) \\ &= (-1)^{(n+1)/2}n \binom{2n}{n} \zeta(n).
    \end{align*}
  On the section $\sum_{i=1}^{2n}x_{i}=1$, this orientation is given by $\mathrm{d}x_{2}\wedge\cdots\wedge\mathrm{d}x_{2n}$. Brown-Schnetz use the projective form
  \[
    \Omega_{W_{n}} = \sum_{i=1}^{2n}(-1)^{i}x_{i}\,\mathrm{d}x_{1}\wedge\cdots\wedge\widehat{\mathrm{d}x_{i}}\wedge\cdots\wedge\mathrm{d}x_{2n}.
  \]
  Substituting $x_{1}=1-\sum_{i=2}^{2n}x_{i}$ gives
  \begin{align*}
    \left.\Omega_{W_{n}}\right|_{\sum_{i=1}^{2n}x_{i}=1} &= -\mathrm{d}x_{2}\wedge\cdots\wedge\mathrm{d}x_{2n},
    \\
    \left.\Omega_{W_{n}}\right|_{x_{2n}=1} &= \mathrm{d}x_{1}\wedge\cdots\wedge\mathrm{d}x_{2n-1}.
  \end{align*}
  Thus, the standard orientation in the chart $x_{2n}=1$ is opposite to the ordered simplex orientation, and the integral with this chart orientation is
  \[
    (-1)^{(n-1)/2}n\binom{2n}{n}\zeta(n).
  \]
  Brown-Schnetz choose the orientation for which the wheel integral is positive \Cite[Introduction and proof of Corollary 63]{brownWheelClassesLocally2025}. Their orientation is therefore $(-1)^{(n-1)/2}$ times the standard orientation in this chart. With their convention, we obtain
  \[
    I_{W_{n}}(\omega^{2n-1}) = n\binom{2n}{n}\zeta(n).
  \]
\end{proof}
\begin{remark}
  Together with \cite[Corollary 8.8]{brownInvariantDifferentialForms2021}, \cref{thm:evaluation-wheel-classes} recovers the known non-vanishing of the odd wheel classes in graph homology: wheels are cycles of degree zero, and $\omega^{2n-1}_{\Lambda_{W_{n}}}$ is a primitive canonical form for odd $n$.
\end{remark}
\section{Gelfand-MacPherson forms}\label{sec:GM-polylogarithms}
We now turn to Gelfand-MacPherson forms. We first recall their construction by integration along Cartan orbits and the resulting differential relation. We then compute the forms associated to powers of the first Pontryagin class in rank two and study their transgression to functions. Finally, we give an integral formula for the forms associated to the Pontryagin character on self-dual Grassmannians.
\subsection{The generic Grassmannian and configuration spaces}
All varieties are defined over $k$ which is either $\mathbf{R}$ or $\mathbf{C}$. Let $\mathrm{G}_{p}^{q} = \mathrm{Gr}(p,p+q)$ denote the Grassmannian of $p$-planes in $(p+q)$-dimensional space. Let $H \subset \mathrm{PGL}_{p+q}$ be the Cartan subgroup given by the image of diagonal matrices in $\mathrm{PGL}_{p+q}$. It is non-canonically isomorphic to $\mathbf{G}_{m}^{p+q-1}$. Denote by $\mathring{\mathrm{G}}_{p}^{q}$ the open of \textsf{generic} $p$-planes, i.e. those $p$-planes that intersect all coordinate $q$-planes only at the origin. In coordinates, $\mathring{\mathrm{G}}_{p}^{q}$ is the locus where all Plücker coordinates do not vanish. The group $H$ preserves $\mathring{\mathrm{G}}_{p}^{q}$ and acts freely on it. Indeed, for a diagonal matrix $\lambda$ with entries $\lambda_{1}, \ldots, \lambda_{p+q}$ the action on $\xi \in \mathrm{G}_{p}^{q}$ is given in Plücker coordinates by  
\[
  \Delta_{I}(\lambda \cdot \xi) = \lambda_{I}\Delta_{I}, \quad \lambda_{I} = \prod_{i\in I}\lambda_{i}
\]
and the claim follows. The quotient may be given two symmetric descriptions, which are exchanged by duality. Let $\xi \in \mathring{\mathrm{G}}_{p}^{q}$ and $e_{1}, \ldots, e_{p+q}$ denote the coordinate basis of $k^{p+q}$. Let $\pi\colon k^{p+q} \to k^{p+q}/\xi$ denote the projection. The genericity condition on $\xi$ implies that every $q$ vectors among the $\pi(e_{i})$ form a basis of $k^{p+q}/\xi$. Identifying $k^{p+q}/\xi$ with $k^{q}$ gives
\[
  H\backslash \mathring{\mathrm{G}}_{p}^{q} \cong \mathrm{Conf}^{\mathrm{gen}}_{p+q}(\mathbf{P}^{q-1})/\mathrm{PGL}_{q}
\]
where generic on the right-hand side means that every $q$ of the $p+q$ points span $\mathbf{P}^{q-1}$. Dually, restrict each $e_{i}^{\vee}$ to $\xi$, obtaining an element of $\xi^{\vee}$. The genericity now says that every $p$ form a basis of $\xi^{\vee}$ and identifying $\xi^{\vee} \cong k^{p}$ gives 
\[
  H\backslash \mathring{\mathrm{G}}_{p}^{q} \cong \mathrm{Conf}^{\mathrm{gen}}_{p+q}(\mathbf{P}^{p-1})/\mathrm{PGL}_{p}
\]
where generic on the right-hand side means that every $p$ of the $p+q$ points span $\mathbf{P}^{p-1}$. We use both incarnations interchangeably, denote the quotient $H\backslash \mathring{\mathrm{G}}_{p}^{q}$ by $\mathrm{C}_{p}^{q}$ and denote by $\pi \colon \mathring{\mathrm{G}}_{p}^{q} \to \mathrm{C}_{p}^{q}$ the quotient map.
\subsection{The bi-Grassmannian complex}
Let $1 \leq j \leq p+q$. Let $\iota_{j}\colon k^{p+q-1} \to k^{p+q}$ denote the inclusion given by setting the $j$th coordinate to be zero. Let $\mathrm{pr}_{j}\colon k^{p+q} \to k^{p+q-1}$ be the map that forgets the $j$th coordinate. For $\xi \in \mathring{\mathrm{G}}_{p}^{q}$ define the maps
\[
  \alpha_{j}(\xi) := \iota_{j}^{-1}(\xi) \in \mathring{\mathrm{G}}_{p-1}^{q}, \quad \beta_{j}(\xi):= \mathrm{pr}_{j}(\xi) \in \mathring{\mathrm{G}}_{p}^{q-1}.
\]
They descend through the $H$-quotients yielding maps (for which we use the same notation)
\[
  \alpha_{j}\colon \mathrm{C}_{p}^{q} \to \mathrm{C}_{p-1}^{q}, \quad \beta_{j}\colon \mathrm{C}_{p}^{q} \to \mathrm{C}_{p}^{q-1}.
\]
If $\mathrm{C}_{p}^{q} \cong \mathrm{Conf}^{\mathrm{gen}}_{p+q}(\mathbf{P}^{q-1})/\mathrm{PGL}_{q}$ then $\alpha_{j}$ forgets the $j$th point and $\beta_{j}$ projects from the $j$th point. Dually, if $\mathrm{C}_{p}^{q} \cong \mathrm{Conf}^{\mathrm{gen}}_{p+q}(\mathbf{P}^{p-1})/\mathrm{PGL}_{p}$ then $\alpha_{j}$ projects from the $j$th point and $\beta_{j}$ forgets the $j$th point. The maps $\alpha_{j}$, $\beta_{j}$ satisfy
\[
  \alpha_{i} \circ \alpha_{j} = \alpha_{j-1} \circ \alpha_{i}, \quad \beta_{i} \circ \beta_{j} = \beta_{j-1} \circ \beta_{i}, \quad i < j
\]
together with the mixed identities
\[
  \alpha_{i} \circ \beta_{j} = \begin{cases}
    \beta_{j-1} \circ \alpha_{i}, &\text{ if }i < j\\
    \beta_{j} \circ \alpha_{i+1}, &\text{ if }i \geq j.
  \end{cases}
\]
\subsection{(Grassmannian) hypersimplices}
The combinatorics of the maps $\alpha_{j}$ and $\beta_{j}$ are controlled by a polytope introduced by Gelfand-MacPherson called the hypersimplex. We denote by $\Delta_{p}^{q}$ the $(p,q)$-\textsf{hypersimplex} defined by
\[
  \Delta_{p}^{q} := \left\{t=(t_{1}, \ldots, t_{p+q}) \in [0,1]^{p+q} \colon \sum_{i=1}^{p+q}t_{i} = p\right\}.
\]
Clearly $\Delta_{1}^{q} = \Delta^{q}$ and $\Delta_{p}^{1} = \Delta^{p}$. There are natural inclusions $a_{j}\colon \Delta_{p-1}^{q} \to \Delta_{p}^{q}$ given by inserting a $1$ on the $j$th coordinate and $b_{j}\colon \Delta_{p}^{q-1} \to \Delta_{p}^{q}$ given by inserting a $0$ on the $j$th coordinate. Namely, the facet given by intersecting $\Delta_{p}^{q}$ with $\{t_{j}=1\}$ is $\Delta_{p-1}^{q}$ and the facet given by intersecting $\Delta_{p}^{q}$ with $\{t_{j}=0\}$ is $\Delta_{p}^{q-1}$. Therefore, $\Delta_{p}^{q}$ is a $(p+q-1)$-dimensional polytope with $2(p+q)$ facets whenever $p,q > 1$. The $(p,q)$-hypersimplex is the moment polytope of the Kähler manifold $\mathrm{G}_{p}^{q}(\mathbf{C})$ with its standard Hamiltonian torus action. In what follows, we write the explicit consequences of \Cite{atiyahConvexityCommutingHamiltonians1982} for this specific case. Indeed, via the Plücker embedding, the \textsf{moment map} $\mu\colon \mathrm{G}_{p}^{q}(\mathbf{C}) \to \mathbf{R}^{p+q}$ is given by 
\begin{equation}\label{eqn:moment-map}
  \mu(\xi) = \frac{\sum_{|I|=p} |\Delta_{I}(\xi)|^{2}e_{I}}{\sum_{|I|=p}|\Delta_{I}(\xi)|^{2}}
\end{equation}
where $e_{I} = \sum_{i\in I}e_{i}$. For an arbitrary $\xi \in \mathrm{G}_{p}^{q}(\mathbf{C})$, the structure of the Cartan-orbits $H(\mathbf{C})\cdot \xi$ is completely described by the \textsf{moment polytope} $P_{\xi} := \mathrm{Conv}\{e_{I}\colon \Delta_{I}(\xi) \neq 0\}$. The Atiyah-Guillemin-Sternberg convexity theorem \Cite{atiyahConvexityCommutingHamiltonians1982,guilleminConvexityPropertiesMoment1982} implies that $\mu(\overline{H(\mathbf{C})\cdot \xi}) = P_{\xi}$. We can say more on the positive real orbit. Let $\xi \in \mathrm{G}_{p}^{q}(\mathbf{R})$. Let $H(\mathbf{R})^{\circ}$ be the connected component of the identify of $H(\mathbf{R})$, i.e. the image of the diagonal matrices with positive entries on $\mathrm{PGL}_{p+q}(\mathbf{R})$. Set $\Gamma_{p}^{q}(\xi) := \overline{H(\mathbf{R})^{\circ} \cdot \xi}$ and $\Gamma_{p}^{q}(\xi)^{\circ} := H(\mathbf{R})^{\circ}\cdot \xi$. Then $\mu$ induces a diffeomorphism \Cite[Theorem 2']{atiyahConvexityCommutingHamiltonians1982} between $\Gamma_{p}^{q}(\xi)^{\circ}$ and $\mathrm{relint}(P_{\xi})$ which extends to a homeomorphism between $\Gamma_{p}^{q}(\xi)$ and $P_{\xi}$. For a generic $\xi$, all Plücker coordinates are non-zero, so that
\[
  P_{\xi} = \mathrm{Conv}\{e_{I}\colon |I|=p\} = \Delta_{p}^{q}.
\]
For non-generic $\xi$, $P_{\xi}$ is a sub-polytope of $\Delta_{p}^{q}$. For $\xi \in \mathring{\mathrm{G}}_{p}^{q}(\mathbf{R})$, we call $\Gamma_{p}^{q}(\xi)$ the \textsf{Grassmannian hypersimplex} of $\xi$. The faces of the Grassmannian hypersimplices are compatible with the bi-Grassmannian complex in an explicit way. Let $a_{j}\colon \mathrm{G}_{p-1}^{q} \hookrightarrow \mathrm{G}_{p}^{q}$ be defined by $a_{j}(\xi) = \mathrm{pr}_{j}^{-1}(\xi)$ and $b_{j}\colon \mathrm{G}_{p}^{q-1} \hookrightarrow \mathrm{G}_{p}^{q}$ given by $b_{j}(\xi) = \iota_{j}(\xi)$. Observe that
\begin{equation}\label{eqn:plucker-compatibility-faces}
  \Delta_{I}(a_{j}(\xi)) = \begin{cases}
    0, &\text{ if } j \notin I\\
    \pm \Delta_{J}(\xi) , &\text{ if } I = J \sqcup \{j\}, 
  \end{cases} \quad \Delta_{I}(b_{j}(\xi)) = \begin{cases}
    0, &\text{ if } j \in I\\
    \Delta_{I}(\xi) , &\text{ if } j \notin I.
  \end{cases}
\end{equation}
From formulas \eqref{eqn:moment-map} and \eqref{eqn:plucker-compatibility-faces} we immediately see that $a_{j} \circ \mu = \mu \circ a_{j}$ and $b_{j} \circ \mu = \mu \circ b_{j}$ which justifies our use of the same notation for them. Taking limits of the Cartan action shows that, for $\xi$ generic
\[
  a_{j}(\Gamma_{p-1}^{q}(\alpha_{j}\xi)) \subset \Gamma_{p}^{q}(\xi), \quad b_{j}(\Gamma_{p}^{q-1}(\beta_{j}\xi)) \subset \Gamma_{p}^{q}(\xi),
\]
so that the diagram
\begin{equation}\label{eqn:facets-grassmannian-hypersimplex}
  \begin{tikzcd}
    \Gamma_{p-1}^{q}(\alpha_{j}\xi) \ar[r, hook, "a_{j}"] \ar[d, "\mu"] & \Gamma_{p}^{q}(\xi) \ar[d, "\mu"] & \Gamma_{p}^{q-1}(\beta_{j}\xi) \ar[l, hook', "b_{j}"'] \ar[d, "\mu"] \\
    \Delta_{p-1}^{q} \ar[r, hook, "a_{j}"] & \Delta_{p}^{q} & \Delta_{p}^{q-1} \ar[l, hook', "b_{j}"']
  \end{tikzcd}
\end{equation}
commutes, i.e. the facets of $\Gamma_{p}^{q}(\xi)$ are exactly $\Gamma_{p-1}^{q}(\alpha_{j}\xi)$ and $\Gamma_{p}^{q-1}(\beta_{j}\xi)$ which recovers \Cite[Theorem 2.3.4]{gelfandGeometryGrassmanniansGeneralization1982}. The Grassmannian hypersimplices admit a simple parametrisation by the hypersimplex. Let $t \in \mathring{\Delta}_{p}^{q}$ and consider the diagonal matrix $\lambda(t)$ with entries $\lambda_{i}(t)= \frac{t_{i}}{1-t_{i}}$. The map $\lambda$ defines a diffeomorphism between $\mathring{\Delta}_{p}^{q}$ and $H(\mathbf{\mathbf{R}})^{\circ}$ \Cite[Proposition 2.3.1]{gelfandGeometryGrassmanniansGeneralization1982}. Consider the map
\begin{equation}\label{eqn:action-hypersimplex}
  \begin{aligned}
  \mathring{\Delta}_{p}^{q} \times \mathring{\mathrm{G}}_{p}^{q}(\mathbf{R}) &\to \mathring{\mathrm{G}}_{p}^{q}(\mathbf{R}) \\
  (t, \xi) &\mapsto t \cdot \xi := \lambda(t)\cdot \xi
  \end{aligned}
\end{equation}
which for fixed $\xi$ induces a diffeomorphism of $\mathring{\Delta}_{p}^{q}$ and $\mathring{\Gamma}_{p}^{q}(\xi)$. Observe that
\begin{equation}\label{eqn:plucker-action}
  \Delta_{I}(t\cdot \xi) = \left(\prod_{i \in I} t_{i}\right)\left(\prod_{i \notin I}1-t_{i}\right)\Delta_{I}(\xi)
\end{equation}
The polynomial formula above shows that the parametrisation \eqref{eqn:action-hypersimplex} extends smoothly to the closed hypersimplex, giving a map $\Delta_{p}^{q} \times \mathring{\mathrm{G}}_{p}^{q} \to \mathrm{G}_{p}^{q}$ which we also denote by $(t,\xi) \mapsto t \cdot \xi$. For fixed $\xi$, we obtain a homeomorphism between $\Delta_{p}^{q}$ and $\Gamma_{p}^{q}(\xi)$, which gives the desired parametrisation of the Grassmannian hypersimplex. Moreover, by looking at Plücker coordinates \eqref{eqn:plucker-action}, the diagram
\[
  \begin{tikzcd}
    \Gamma_{p-1}^{q}(\alpha_{j}\xi) \ar[r, hook, "a_{j}"]  & \Gamma_{p}^{q}(\xi)  & \Gamma_{p}^{q-1}(\beta_{j}\xi) \ar[l, hook', "b_{j}"']  \\
    \Delta_{p-1}^{q} \ar[u] \ar[r, hook, "a_{j}"] & \Delta_{p}^{q} \ar[u] & \Delta_{p}^{q-1} \ar[l, hook', "b_{j}"'] \ar[u]
  \end{tikzcd}
\]
commutes.
\subsection{Gelfand-MacPherson forms}
Let $P$ be a stable characteristic class on real vector bundles with values in $\mathrm{H}_{\mathrm{sing}}^{\ell}(-;\mathbf{R})$, i.e. $P$ is a polynomial on Pontryagin classes. Let $\gamma_{p}^{q} \to \mathrm{G}_{p}^{q}(\mathbf{R})$ be the tautological bundle. Let $P_{p}^{q} \in \Omega^{\ell}(\mathrm{G}_{p}^{q}(\mathbf{R}))^{\mathrm{SO}(p+q)}$ be the unique $\mathrm{SO}(p+q)$-invariant differential form which represents $P(\gamma_{p}^{q})$ in de Rham cohomology (which exists because $\mathrm{G}_{p}^{q}(\mathbf{R})$ is a symmetric space, it is the harmonic representative for an $\mathrm{SO}(p+q)$-invariant metric). A fortiori $P_{p}^{q}$ happens to also be $\mathrm{O}(p+q)$-invariant. Indeed, let $g \in \mathrm{O}(p+q)$, then $g^{*}P_{p}^{q}$ is also $\mathrm{SO}(p+q)$-invariant because the group $\mathrm{SO}(p+q)$ is normal in $\mathrm{O}(p+q)$. But $g^{*}\gamma_{p}^{q} = \gamma_{p}^{q}$, so $g^{*}P_{p}^{q}$ is another $\mathrm{SO}(p+q)$-invariant representative of $P(\gamma_{p}^{q})$, so by uniqueness $P_{p}^{q} = g^{*}P_{p}^{q}$.
\begin{definition}\cite[Definition 1.2.1]{gelfandGeometryGrassmanniansGeneralization1982}\label{def:Gelfand-MacPherson-polylogarithmic-forms}
  The \textsf{Gelfand-MacPherson form} $\widehat{P}_{p}^{q} \in \Omega^{\ell-(p+q-1)}(\mathrm{C}_{p}^{q}(\mathbf{R}))$ is the integration along the fibres $\pi_{*}(P_{p}^{q}\lvert_{\mathring{\mathrm{G}}_{p}^{q}})$ of $\pi\colon \mathring{\mathrm{G}}_{p}^{q}(\mathbf{R}) \to \mathrm{C}_{p}^{q}(\mathbf{R})$.
\end{definition}
The previous definition is deserving of an explanation, since clearly the fibres of $\pi$ are non-compact (Gelfand-MacPherson circumvent this via ``enhanced'' configurations \Cite[Section 2.4.1]{gelfandGeometryGrassmanniansGeneralization1982} which we rather not introduce). The simple explanation is that over every configuration of points, the fibres of $\pi$ admit a canonical compactification, i.e. the Grassmannian hypersimplex. First, observe that the principal $H(\mathbf{R})$-bundle $\pi\colon \mathring{\mathrm{G}}_{p}^{q} \to \mathrm{C}_{p}^{q}(\mathbf{R})$ is trivial. Write $\mathrm{C}_{p}^{q}(\mathbf{R})$ as $\mathrm{Conf}^{\mathrm{gen}}_{p+q}(\mathbf{RP}^{q-1})/\mathrm{PGL}_{q}(\mathbf{R})$. The first $q+1$ points of every generic configuration of $p+q$ points in $\mathbf{RP}^{q-1}$ form a projective frame. The space of ordered projective frames on $\mathbf{RP}^{q-1}$ is a smooth $\mathrm{PGL}_{q}(\mathbf{R})$-torsor. Hence, every configuration admits a unique representative whose first $q+1$ points are a fixed standard frame. The remaining points have non-zero homogeneous coordinates with respect to the fixed standard frame. By normalising, we may choose representatives for all the $p+q$ points and any choice of normalisation lies on the same $H(\mathbf{R})$-orbit, which shows that $\pi$ admits a global smooth section. Set $F\subset H(\mathbf{R})$ for the image of the diagonal sign matrices, so that $F \cong (\mathbf{Z}/2 \mathbf{Z})^{p+q-1}$ non-canonically. We can then write $H(\mathbf{R}) = H(\mathbf{R})^{\circ} \times F$. By triviality of $\pi$, choose a smooth global section $s\colon \mathrm{C}_{p}^{q}(\mathbf{R}) \to \mathring{\mathrm{G}}_{p}^{q}(\mathbf{R})$. From \eqref{eqn:plucker-action}, there is a diffeomorphism
\[
  \begin{aligned}
    F \times \mathring{\Delta}_{p}^{q} \times \mathrm{C}_{p}^{q}(\mathbf{R}) &\xrightarrow{\cong} \mathring{\mathrm{G}}_{p}^{q}(\mathbf{R}) \\ (\varepsilon, t, c) &\mapsto t \cdot \varepsilon s(c).
  \end{aligned}
\]
For $\varepsilon \in F$, let $r_{\varepsilon} \colon \Delta_{p}^{q} \times \mathrm{C}_{p}^{q}(\mathbf{R}) \to \mathrm{G}_{p}^{q}(\mathbf{R})$ be the unique extension of $(t, c) \mapsto t \cdot \varepsilon s(c)$. The image of $r_{\varepsilon}$ over a fixed $c \in \mathrm{C}_{p}^{q}(\mathbf{R})$ is $\Gamma_{p}^{q}(\varepsilon s(c))$. Let $q\colon \Delta_{p}^{q} \times \mathrm{C}_{p}^{q}(\mathbf{R}) \to \mathrm{C}_{p}^{q}(\mathbf{R})$ be the projection onto the second factor, which clearly has compact fibres. We orient the fibres by the images of $E_{11},\ldots,E_{p+q-1,p+q-1}$ in the diagonal Lie algebra modulo scalar matrices, as in \Cite[Section 1.1.2]{gelfandGeometryGrassmanniansGeneralization1982}, and use fibre-first integration. Here $E_{ij}$ denotes the $(p+q)\times(p+q)$ matrix with a $1$ in position $(i,j)$ and zeros elsewhere. Here $q_{*}$ uses the hypersimplex orientation fixed below, which differs from the Cartan fibre orientation by $(-1)^{p+q-1}$. We define
\[
  \widehat{P}_{p}^{q} = \pi_{*}\left(P_{p}^{q}\lvert_{\mathring{\mathrm{G}}_{p}^{q}}\right) := (-1)^{p+q-1}\sum_{\varepsilon \in F} q_{*}r_{\varepsilon}^{*}P_{p}^{q}.
\]
However, since $P_{p}^{q}$ is $\mathrm{O}(p+q)$-invariant and $F \subset \mathrm{O}(p+q)$, plus the fact that $H(\mathbf{R})$ is abelian, so sign matrices preserve the fibre orientation, we obtain that
\begin{equation}\label{eqn:integration-over-fibre}
  (\widehat{P}_{p}^{q})_{c} = (-1)^{p+q-1}2^{p+q-1}q_{*}r_{1}^{*}P_{p}^{q}.
\end{equation}
Any other choice of section $s'$ gives the same family of Grassmannian hypersimplices up to $\mathrm{O}(p+q)$, so that $\widehat{P}_{p}^{q}$ is independent of the choice of section.
\subsection{The differential relation}
Once and for all we fix an orientation on the hypersimplex $\Delta_{p}^{q}$ by the ordered affine coordinates $t_{2}-t_{1}, t_{3}-t_{2}, \ldots, t_{p+q}-t_{p+q-1}$. With this convention, the boundary orientation of the facet $\{t_{j}=1\}$ differs from the fixed orientation of $\Delta_{p-1}^{q}$ by $(-1)^{j}$ while for $\{t_{j}=0\}$ the difference to the orientation on $\Delta_{p}^{q-1}$ is $-(-1)^{j}$. Our fixed orientation defines a generator of $\mathrm{H}^{\mathrm{sing}}_{p+q-1}(\Delta_{p}^{q}, \partial \Delta_{p}^{q};\mathbf{Z})$ and it has homological boundary \Cite[Proposition 2.1.3]{gelfandGeometryGrassmanniansGeneralization1982}
\begin{equation}\label{eqn:homological-boundary}
  \partial [\Delta_{p}^{q}] = \sum_{j=1}^{p+q}(-1)^{j}a_{j*}[\Delta_{p-1}^{q}] - \sum_{j=1}^{p+q}(-1)^{j}b_{j*}[\Delta_{p}^{q-1}].
\end{equation}
From the commutativity of \eqref{eqn:facets-grassmannian-hypersimplex} and \eqref{eqn:homological-boundary} we see that for $\xi$ generic, orienting $\Gamma_{p}^{q}(\xi)$ by the Cartan fibre orientation,
\begin{equation}\label{eqn:homological-boundary-grassmannian}
  \partial [\Gamma_{p}^{q}(\xi)] = -\sum_{j=1}^{p+q}(-1)^{j}a_{j*}[\Gamma_{p-1}^{q}(\alpha_{j}\xi)] + \sum_{j=1}^{p+q}(-1)^{j}b_{j*}[\Gamma_{p}^{q-1}(\beta_{j}\xi)].
\end{equation}
Since $P$ is a stable characteristic class, universality gives $a_{j}^{*}P_{p}^{q} = P_{p-1}^{q}$ and $b_{j}^{*}P_{p}^{q} = P_{p}^{q-1}$. We triangulate the hypersimplex and apply Stokes' theorem to each simplex; the contributions from shared faces cancel. Thus \eqref{eqn:homological-boundary-grassmannian} gives for $p>1$ and $q>1$ the differential relation \Cite[Theorem 1.2.2]{gelfandGeometryGrassmanniansGeneralization1982}\footnote{With the fibre orientation of \Cite[Section 1.1.2]{gelfandGeometryGrassmanniansGeneralization1982} and fibre-first integration, the relation has an additional factor $(-1)^{p+q-1}$ compared with their Theorem 1.2.2. Fibrewise Stokes contributes $(-1)^{p+q}$, while comparison with the hypersimplex orientations reverses the facet signs. The resulting functional equations are unchanged.}
\begin{equation}\label{eqn:differential-relation}
  \mathrm{d}\widehat{P}_{p}^{q} = 2(-1)^{p+q-1}\left(\sum_{j=1}^{p+q}(-1)^{j}\alpha_{j}^{*}\widehat{P}_{p-1}^{q} - \sum_{j=1}^{p+q}(-1)^{j}\beta_{j}^{*}\widehat{P}_{p}^{q-1}\right)
\end{equation}
where the mysterious $2$ is just $\frac{2^{p+q-1}}{2^{p+q-2}}$ in disguise.
\subsection{Vanishing and functional equations}
The topological problem is the following: fix a given stable class $P$ with values in $\mathrm{H}^{\ell}_{\mathrm{sing}}(-;\mathbf{R})$. Find a pair $(p,q)$ such that $\widehat{P}_{p}^{q}$ is closed, so that the left-hand side of the differential relation \eqref{eqn:differential-relation} is zero, giving rise to a functional equation on $\widehat{P}_{p-1}^{q}$ and $\widehat{P}_{p}^{q-1}$. 

Surprisingly, the form $P_{p}^{q}$ always vanishes in $H(\mathbf{R})$-orbits. In \Cite[Theorem 1.2.3]{gelfandGeometryGrassmanniansGeneralization1982} they show this via direct calculation, but as explained in \cite{atiyahConvexityCommutingHamiltonians1982} this may be seen as a consequence of the Kähler structure in $\mathrm{G}_{p}^{q}(\mathbf{C})$. Indeed, since $\mathrm{G}_{p}^{q}(\mathbf{C})$ is Kähler and its torus action admits fixed points, then \Cite[Theorem 4]{atiyahConvexityCommutingHamiltonians1982} holds and all harmonic forms of positive degree vanish on $H(\mathbf{R})$-orbits. In particular
\begin{theorem}\label{thm:vanishing-orbits}
For every homogeneous polynomial $P$ of positive degree in Pontryagin classes, $P_{p}^{q}$ vanishes on every $H(\mathbf{R})$-orbit.
\end{theorem}
\begin{corollary}\label{cor:0-form-vanishing}
  Whenever $\ell = p+q-1$, $\widehat{P}_{p}^{q}$ vanishes. In particular
  \[
    0 = \sum_{j=1}^{p+q}(-1)^{j}\alpha_{j}^{*}\widehat{P}_{p-1}^{q} - \sum_{j=1}^{p+q}(-1)^{j}\beta_{j}^{*}\widehat{P}_{p}^{q-1}.
  \]
\end{corollary}
\subsection{The first Pontryagin class and Rogers dilogarithm}
For the purposes of exposition, we carry out in detail the consequences of \cref{cor:0-form-vanishing} in the first non-trivial case. This subsection essentially follows \Cite[Section 0]{gelfandGeometryGrassmanniansGeneralization1982}. Set $P = p_{1}$, the first Pontryagin class. Then $\ell = \deg P = 4$. First, observe that $\mathrm{G}_{1}^{3}(\mathbf{R}) \cong \mathbf{RP}^{3}$, so its tautological bundle has vanishing Pontryagin classes. In particular, $\widehat{P}_{1}^{3} = 0$. We obtain the following functional equation for $\widehat{P}_{2}^{2}$ on $\mathrm{C}_{2}^{3}(\mathbf{R})$ via \cref{cor:0-form-vanishing}
\begin{equation}\label{eqn:functional-eqn-G_2^2}
  0 = \sum_{j=1}^{5}(-1)^{j}\beta_{j}^{*}\widehat{P}_{2}^{2}.
\end{equation}
Recall that we may identify $\mathrm{C}_{2}^{q}$ with $\mathrm{Conf}^{\mathrm{gen}}_{2+q}(\mathbf{P}^{1})/\mathrm{PGL}_{2}$. However, in this specific case $\mathrm{Conf}^{\mathrm{gen}}_{2+q}(\mathbf{P}^{1})/\mathrm{PGL}_{2} = \mathcal{M}_{0,2+q}$, the moduli space of genus $0$ curves with $(2+q)$ marked points. The maps $\beta_{j}\colon \mathcal{M}_{0,2+q} \to \mathcal{M}_{0,1+q}$ correspond to forgetting the $j$th marked point. In coordinates, normalise $(t_{1}, \ldots, t_{5}) = (\infty, 0, 1, x, y)$, so that $\mathcal{M}_{0,5}(\mathbf{R}) \cong \{(x,y) \in (\mathbf{R}\setminus\{0,1\})^{2} \colon x \neq y\}$, $\mathcal{M}_{0,4}(\mathbf{R}) \cong \mathbf{R}\setminus \{0,1\}$ and the maps $ \beta_{j}\colon \mathcal{M}_{0,5}(\mathbf{R}) \to \mathcal{M}_{0,4}(\mathbf{R})$ are
\[
  \beta_{1}(x,y) = \frac{x(y-1)}{y(x-1)}, \quad \beta_{2}(x,y)=\frac{y-1}{x-1}, \quad \beta_{3}(x,y)=\frac{y}{x}, \quad \beta_{4}(x,y)=y, \quad \beta_{5}(x,y)=x.
\]
Since $\dim \mathcal{M}_{0,4}(\mathbf{R}) = 1$, then $\widehat{P}_{2}^{2}$ is closed. In every connected component of $\mathcal{M}_{0,4}(\mathbf{R})$ we may choose a primitive $\Phi$, so that $\mathrm{d}\Phi = \widehat{P}_{2}^{2}$. Therefore, from \eqref{eqn:functional-eqn-G_2^2}, we obtain
\begin{equation}\label{eqn:5-term-relation}
  \mathrm{constant} = -\Phi\left(\frac{x(y-1)}{y(x-1)}\right) + \Phi\left(\frac{y-1}{x-1}\right) - \Phi\left(\frac{y}{x}\right) + \Phi(y) - \Phi(x).
\end{equation}
which is known as the $5$-term relation. Consider the real Rogers dilogarithm defined in $\mathbf{R}\setminus \{0,1\}$
\begin{equation}\label{eqn:rogers-dilog}
  L_{\mathbf{R}}(x) := \mathrm{Li}_{2}^{\mathbf{R}}(x) + \frac{1}{2}\log|x|\log|1-x|, \quad \mathrm{Li}_{2}^{\mathbf{R}}(x) := - \int_{0}^{x} \frac{\log|1-t|}{t}\mathrm{d}t,
\end{equation}
where $\mathrm{Li}_{2}^{\mathbf{R}}(x) = \mathfrak{R}\mathrm{Li}_{2}(x)$ is the real part of the classical dilogarithm function. In \cref{sec:rank-2-poly-forms}, we compute by direct integration that $\Phi=\frac{2}{\pi^{2}}L_{\mathbf{R}}+C$, where $C$ is constant on each connected component.
\subsection{Powers of the first Pontryagin class and Damiano vanishing}
We now work on $\mathrm{G}_{2}^{q}$. On the oriented double cover of $\mathrm{G}_{2}^{q}(\mathbf{R})$, consider $e$ the invariant Euler $2$-form of the tautological oriented $2$-plane bundle. Universally, we have that $p_{1}=e^{2}$. Let $V$ be the tangent distribution to the $H(\mathbf{R})$-orbits. We may decompose $e$, relative to any complement of $V$, according to the vertical degree $e = e_{2,0} + e_{1,1} + e_{0,2}$, where $e_{2,0}$ is the purely vertical component of $e$. By \cref{thm:vanishing-orbits}, $p_{1}(\Lambda^{4}V) = 0$. Since $p_{1}=e^{2}$, then $e^{2}_{2,0} = 0$. Now we look at $p_{1}^{k} = e^{2k}$. A monomial in the expansion of $e^{2k}$ which contains at most one $e_{2,0}$ has vertical degree at most $2 + (2k-1) = 2k+1$. Any term of vertical degree at least $2k+2$ must therefore contain at least two copies of $e_{2,0}$, so it vanishes. In conclusion,
\[
  (e^{2k})_{\text{vertical degree }\geq 2k+2} = 0.
\]
However, the fibre of $\pi\colon \mathring{\mathrm{G}}_{2}^{q}(\mathbf{R}) \to \mathcal{M}_{0,2+q}(\mathbf{R})$ has dimension $q+1$, so after fibre integration we obtain a vanishing result due to D.B. Damiano \cite[Theorem 5.1]{damianoWebsCharacteristicForms1983}\Cite[Theorem 6.3.1]{damianoWEBSABELIANEQUATIONS1980} which he obtained by direct computation of the Euler class.
\begin{theorem}\label{thm:damiano-vanishing}
  Let $P=p_{1}^{k}$ be the $k$th power of the first Pontryagin class. Then $\widehat{P}_{2}^{2k+1}=0$. 
\end{theorem}
Since the tautological bundle of $\mathrm{G}_{1}^{2k+1}(\mathbf{R})\cong \mathbf{RP}^{2k+1}$ has no Pontryagin classes, we obtain the following functional equation from \eqref{eqn:differential-relation}.
\begin{corollary}\label{cor:functional-equation-l-power-first-pontryagin-class}
  Let $P=p_{1}^{k}$ be the $k$th power of the first Pontryagin class. Then the Gelfand-MacPherson form $\widehat{P}_{2}^{2k}$ satisfies the functional equation
  \[
    0 = \sum_{j=1}^{2k+3}(-1)^{j}\beta_{j}^{*}\widehat{P}_{2}^{2k}.
  \]
\end{corollary}
\subsection{Gelfand-MacPherson forms in rank $2$}\label{sec:rank-2-poly-forms}
In this subsection we compute explicitly the Gelfand-MacPherson form $\widehat{P}_{2}^{2k}$ for $P=p_{1}^{k}$. First, we need an expression for the $k$th power $p_{1}^{k}$ of the first Pontryagin class of the tautological bundle of $\mathrm{G}_{2}^{2k}(\mathbf{R})$. The special feature of this particular case is that $p_{1}^{k}$ is an invariant volume form for $\mathrm{G}_{2}^{2k}(\mathbf{R})$. Consider the graph chart about the standard plane $\mathbf{R}^{2} \subset \mathbf{R}^{2} \oplus \mathbf{R}^{2k}$, i.e.
\begin{equation}\label{eqn:graph-coordinates}
  X = (x_{ij}) \in \mathrm{Mat}_{2\times 2k}(\mathbf{R}) \mapsto \mathrm{rowspan}(\mathrm{I}_{2\times 2}\mid X)
\end{equation}
oriented by
\[
  \mathrm{d}X := \bigwedge_{j=1}^{2k}\mathrm{d}x_{1j} \wedge \mathrm{d}x_{2j}.
\]
\begin{proposition}\label{prop:first-pontryagin-class}
  Let $P=p_{1}^{k}$ be the $k$th power of the first Pontryagin class. In the graph coordinates \eqref{eqn:graph-coordinates}, the $\mathrm{SO}(2+2k)$-invariant expression for $P_{2}^{2k}$ is
  \[
    P_{2}^{2k} = \frac{(2k)!}{(2\pi)^{2k}}\frac{\mathrm{d}X}{\mathrm{det}(\mathrm{I}_{2\times 2} + XX^{t})^{k+1}}.
  \]
\end{proposition}
\begin{proof}
  We first determine an $\mathrm{SO}(2+2k)$-invariant volume form on $\mathrm{G}_{2}^{2k}(\mathbf{R})$. We may carry out the calculation first at a point, and then use the homogeneous space structure. At the standard plane $o$ (which in graph coordinates corresponds to $X = 0$), identify
  \[
    \mathrm{T}_{o}\mathrm{G}_{2}^{2k}(\mathbf{R}) \cong \mathrm{Hom}(\mathbf{R}^{2}, \mathbf{R}^{2k}) \cong \mathrm{Mat}_{2\times 2k}(\mathbf{R})
  \]
  with invariant metric $g_{o}(U,V) = \mathrm{tr}(UV^{t})$. In graph coordinates, we obtain
  \[
    g_{X}(U,V) = \mathrm{tr}( (\mathrm{I}_{2\times 2} + XX^{t})^{-1}U(\mathrm{I}_{2k \times 2k}+ X^{t}X)^{-1}V^{t}).
  \]
  Observe that
  \[
    \det g_{X} = \det(\mathrm{I}_{2k \times 2k}+ X^{t}X)^{-2} \det( \mathrm{I}_{2\times 2} + XX^{t} )^{-2k}.
  \]
  However, by Sylvester's determinant identity,
  \[
    \det(\mathrm{I}_{2k \times 2k}+ X^{t}X) = \det( \mathrm{I}_{2\times 2} + XX^{t} ),
  \]
  so that
  \[
    \mathrm{d}\mathrm{vol}_{g} = \sqrt{\det g_{X}} \mathrm{d}X = \frac{\mathrm{d}X}{\mathrm{det}(\mathrm{I}_{2\times 2} + XX^{t})^{k+1}}
  \]
  is an $\mathrm{SO}(2+2k)$-volume form for $\mathrm{G}_{2}^{2k}(\mathbf{R})$. To identify the constant, we proceed via Chern-Weil as in \Cite[Proof of Theorem 5.1.]{damianoWebsCharacteristicForms1983}. Pass to the oriented double cover of the Grassmannian, so that the tautological $2$-bundle is oriented. Let $\Omega$ be the curvature of its canonical connection. In an oriented orthonormal frame we may write
  \[
    \Omega = \begin{pmatrix}
      0 & \Omega_{12} \\
      - \Omega_{12} & 0
    \end{pmatrix}
  \]
  so that $e = \frac{1}{2\pi}\Omega_{12}$. Since $p_{1}= e^{2}$, then
  \[
    P_{2}^{2k} = \frac{1}{(2\pi)^{2k}}\Omega_{12}^{2k}.
  \]
  Again, we may simply evaluate at $o$. Let $\theta_{ij}$ be the orthonormal coframe of $\mathrm{Hom}(\mathbf{R}^{2}, \mathbf{R}^{2k})$, so that at the origin $X=0$ we have
  \[
    \Omega_{12} = - \sum_{j} \theta_{1j}\wedge \theta_{2j}
  \]
  because $\Omega = -\theta \wedge \theta^{t}$. Therefore, at $X=0$,
  \[
    P_{2}^{2k}\lvert_{X=0} = \frac{(2k)!}{(2\pi)^{2k}}\mathrm{d}X
  \]
  and the claim follows using the homogeneous structure.
\end{proof}
In order to state our main result, we need to discuss some background. Write $\mathrm{C}_{2}^{2k}(\mathbf{R}) = \mathcal{M}_{0,2+2k}(\mathbf{R})$. Consider configurations of the form $x = (\infty, x_{1}, \ldots, x_{2k+1})$ with $x_{i} \in \mathbf{R}$ and $x_{i} \neq x_{j}$. In other words, we make the identification
\[
  \mathcal{M}_{0,2+2k}(\mathbf{R}) \cong \mathrm{Conf}_{2k+1}(\mathbf{R})/\mathrm{Aff}_{1}(\mathbf{R}).
\]
Choose homogeneous representatives
\[
  v_{0} = \begin{pmatrix}
    1 \\ 0
  \end{pmatrix}, \quad v_{i} = \begin{pmatrix}
    x_{i} \\ 1
  \end{pmatrix}, \quad 1 \leq i \leq 2k+1
\]
so that the $2 \times (2 + 2k)$ matrix $S = (v_{0}, \ldots, v_{2k+1})$ has image $x$ under $\pi \colon \mathring{\mathrm{G}}_{2}^{2k}(\mathbf{R}) \to \mathrm{C}_{2}^{2k}(\mathbf{R})$. Its Plücker coordinates $\Delta_{ij} = \det(v_{i}, v_{j})$ are
\[
  \Delta_{0i}=1, \quad \Delta_{ij} = x_{i}-x_{j}\quad (i,j > 0).
\]
The action of $H(\mathbf{R})$ re-scales the columns independently. Therefore, for $r_{i} > 0$, set
\[
  S(r) = (\sqrt{r_{0}}v_{0}, \ldots, \sqrt{r_{2k+1}}v_{2k+1})
\]
which parametrises the $H(\mathbf{R})^{\circ}$-orbit of $S$. For a $2 \times (2+2k)$ matrix $A$ its Gram determinant is $Q(A) = \det(AA^{t})$ where $G(A) = AA^{t}$ is its Gram matrix. By Cauchy-Binet,
\begin{equation}\label{eqn:Gram-determinant-rank-2}
  Q(S(r)) = \det(S(r)S(r)^{t}) = \sum_{i<j}r_{i}r_{j}\Delta_{ij}^{2} = r_{0}\sum_{i=1}^{2k+1}r_{i} + \sum_{1 \leq i < j \leq 2k+1} r_{i}r_{j}(x_{i}-x_{j})^{2}.
\end{equation}
We are now ready to give a completely explicit formula for $\widehat{P}_{2}^{2k}$ for $P=p_{1}^{k}$ the $k$th power of the first Pontryagin class.
\begin{theorem}\label{thm:Gelfand-MacPherson-polylogarithmic-form-rank-2}
  Let $P=p_{1}^{k}$ be the $k$th power of the first Pontryagin class. For $x = (x_{1},\ldots,x_{2k+1}) \in \mathrm{Conf}_{2k+1}(\mathbf{R})/\mathrm{Aff}_{1}(\mathbf{R}) \cong \mathrm{C}_{2}^{2k}(\mathbf{R})$, the Gelfand-MacPherson form $\widehat{P}_{2}^{2k} \in \Omega^{2k-1}(\mathrm{C}_{2}^{2k}(\mathbf{R}))$ is
  \begin{equation}\label{eqn:affine-expression-rank-2}
     \widehat{P}_{2}^{2k} = -\frac{4}{(2\pi)^{2k}}\sum_{i<j}(-1)^{i+j}\log|x_{i}-x_{j}|\bigwedge_{\ell\neq i,j}\mathrm{d}\log \left|\frac{x_{\ell} - x_{i}}{x_{\ell} - x_{j}}\right|
  \end{equation}
  Setting $x_{1} = 0$ and $x_{2}=1$ to fix an $\mathrm{Aff}_{1}(\mathbf{R})$-gauge, i.e. so that $x_{3}, \ldots, x_{2k+1}$ are coordinates on $\mathrm{C}_{2}^{2k}(\mathbf{R})$, gives
  \begin{equation}\label{eqn:gauge-fixed-rank-2}
    \widehat{P}_{2}^{2k} = -\frac{4}{(2\pi)^{2k}}\sum_{i<j}\frac{(x_{i}-x_{j})^{2k-2}\log|x_{i}-x_{j}|}{\prod_{\ell \neq i,j}(x_{i}-x_{\ell})(x_{j}-x_{\ell})}\mathrm{d}x_{3} \wedge \cdots \wedge \mathrm{d}x_{2k+1}.
  \end{equation}
  
\end{theorem}
\begin{proof}
  We compute on the affine slice $x_{1}=0$, $x_{2}=1$ and arrive at \eqref{eqn:gauge-fixed-rank-2} first. Points of the positive Cartan fibre over this configuration may be represented by
  \[
    S(a) = \begin{pmatrix}
      a_{0} & 0 & 1 & a_{3}x_{3} & \cdots & a_{2k+1}x_{2k+1} \\
      0 & a_{1} & 1 & a_{3} & \cdots & a_{2k+1}
    \end{pmatrix}, \quad a_{i} > 0
  \]
  and $a_{2}=1$. To change to graph coordinates $X$ \eqref{eqn:graph-coordinates} we multiply on the left by $\mathrm{diag}(a_{0}^{-1}, a_{1}^{-1})$. A direct calculation of the Jacobian of this change of coordinates gives
  \[
    \frac{\mathrm{d}X}{\det(\mathrm{I}_{2\times 2} + XX^{t})^{k+1}} = (-1)^{k}\frac{a_{0}a_{1}a_{3}\cdots a_{2k+1}}{Q(S(a))^{k+1}}\mathrm{d}_{\mathrm{fib}}a \wedge \mathrm{d}x_{3} \wedge \cdots \wedge \mathrm{d}x_{2k+1}.
  \]
  where $\mathrm{d}_{\mathrm{fib}}a$ is the coordinate volume form in coordinates $a_{0}, a_{1}, a_{3}, \ldots, a_{2k+1}$ and $Q$ is the Gram determinant \eqref{eqn:Gram-determinant-rank-2}. Set $r_{i}=a_{i}^{2}$ for all $0\leq i \leq 2k+1$. Then
  \[
    a_{0}a_{1}a_{3} \cdots a_{2k+1} \mathrm{d}_{\mathrm{fib}}a = \frac{1}{2^{2k+1}}\mathrm{d}_{\mathrm{fib}}r.
  \]
  These coordinates have the opposite orientation to the Cartan fibre orientation; the integrals below use the positive coordinate measures, so this contributes an additional minus sign. From \cref{prop:first-pontryagin-class} and \eqref{eqn:integration-over-fibre}, the power $2^{2k+1}$ cancels and we obtain
  \begin{equation}\label{eqn:GM-rank-2-r-variables}
    \widehat{P}_{2}^{2k} = \frac{(-1)^{k+1}(2k)!}{(2\pi)^{2k}}\left(\int_{\mathbf{R}^{2k+1}_{>0}} \frac{\mathrm{d}_{\mathrm{fib}}r}{Q(S(r))^{k+1}}\right) \mathrm{d}x_{3} \wedge \cdots \wedge  \mathrm{d}x_{2k+1}.
  \end{equation}
   Write $Q(S(r))$ as $r_{0}L + V$ where
  \[
    L:= \sum_{i=1}^{2k+1} r_{i}, \quad V:= \sum_{1 \leq i < j \leq 2k+1} r_{i}r_{j}(x_{i}-x_{j})^{2}.
  \]
  The $r_{0}$-integration is elementary
  \[
    \int_{0}^{\infty} \frac{\mathrm{d}r_{0}}{(r_{0}L + V)^{k+1}} = \frac{1}{kLV^{k}},
  \]
  which leaves us with
  \[
    \frac{(-1)^{k+1}(2k)!}{k(2\pi)^{2k}} \int_{\mathbf{R}^{2k}_{>0}} \frac{\mathrm{d}r_{1} \mathrm{d}r_{3} \cdots \mathrm{d}r_{2k+1}}{\left(\sum_{i=1}^{2k+1}r_{i}\right)\left(\sum_{i<j}r_{i}r_{j}(x_{i}-x_{j})^{2}\right)^{k}}.
  \]
  The remaining integral is projective in the variables $r_{1}, \ldots, r_{2k+1}$. Set
  \[
    t_{i}:= \frac{r_{i}}{L}, \quad \sum_{i=1}^{2k+1}t_{i}=1,
  \]
  so that the $t_{i}$ are coordinates in the simplex $\Delta^{2k}$. Since $r_{2}=1$, we have that $L=t_{2}^{-1}$ and
  \[
    \mathrm{d}r_{1}\mathrm{d}r_{3}\cdots \mathrm{d}r_{2k+1} = L^{2k+1} \mathrm{d}t.
  \]
  Observe that
  \[
    L\left(\sum_{i<j}r_{i}r_{j}(x_{i}-x_{j})^{2}\right)^{k} = L^{2k+1}\left(\sum_{i<j}t_{i}t_{j}(x_{i}-x_{j})^{2}\right)^{k},
  \]
  so that the powers of $L$ cancel. We arrive at the following form of \eqref{eqn:GM-rank-2-r-variables}
  \[
    \widehat{P}_{2}^{2k} = \frac{(-1)^{k+1}(2k)!}{k(2\pi)^{2k}} \left(\int_{\Delta^{2k}} \frac{\mathrm{d}t}{\left(\sum_{i<j}t_{i}t_{j}(x_{i}-x_{j})^{2}\right)^{k}}\right)\mathrm{d}x_{3} \wedge \cdots \wedge \mathrm{d}x_{2k+1}.
  \]
  The integral in parentheses is exactly \cref{thm:quadratic-barycentric-transgression}, so that
  \[
    \widehat{P}_{2}^{2k} = -\frac{4}{(2\pi)^{2k}}\sum_{i<j}\frac{(x_{i}-x_{j})^{2k-2}\log|x_{i}-x_{j}|}{\prod_{\ell \neq i,j}(x_{i}-x_{\ell})(x_{j}-x_{\ell})}\mathrm{d}x_{3} \wedge \cdots \wedge \mathrm{d}x_{2k+1}
  \]
  as desired. Now we show how \eqref{eqn:gauge-fixed-rank-2} gives \eqref{eqn:affine-expression-rank-2}. On $\mathrm{Conf}_{2k+1}(\mathbf{R})$ set
  \[
    \Omega = \mathrm{d}x_{1} \wedge \cdots \wedge \mathrm{d}x_{2k+1}, \quad T=\sum_{i}\partial_{i}, \quad D=\sum_{i}x_{i}\partial_{i}, \quad \nu = \iota_{D}\iota_{T}\Omega
  \]
  where $T$ and $D$ are the infinitesimal generators (of translation and dilation) of the action of $\mathrm{Aff}_{1}(\mathbf{R})$ on $\mathrm{Conf}_{2k+1}(\mathbf{R})$. It is enough to show that the following identity holds
  \[
    \frac{(x_{i}-x_{j})^{2k-2}}{\prod_{\ell\neq i,j} (x_{i}-x_{\ell})(x_{j}-x_{\ell})} \nu = (-1)^{i+j}\bigwedge_{\ell\neq i,j} \mathrm{d}\log \left|\frac{x_{\ell} - x_{i}}{x_{\ell}-x_{j}}\right|. 
  \]
  Indeed, both sides are invariant under $\mathrm{Aff}_{1}(\mathbf{R})$, and a direct check shows that they agree on the gauge $x_{i} = 0$ and $x_{j}=1$. Putting $x_{1}=0$ and $x_{2}=1$ gives $\nu = \mathrm{d}x_{3} \wedge \cdots \wedge \mathrm{d}x_{2k+1}$ and \eqref{eqn:affine-expression-rank-2} follows.
\end{proof}
\begin{remark}\label{rem:specialisation-rogers-dilog}
  Set $k=1$ in \eqref{eqn:gauge-fixed-rank-2} with $x:=x_{3}$. We obtain (see \eqref{eqn:rogers-dilog})
  \[
    \widehat{P}_{2}^{2} = - \frac{1}{\pi^{2}}(-\log|x|\mathrm{d}\log|1-x| + \log|1-x|\mathrm{d}\log|x|) = \frac{2}{\pi^{2}} \mathrm{d}L_{\mathbf{R}}.
  \]
\end{remark}
\begin{remark}\label{rem:symmetric-expression-rank-2}
  In the form of \eqref{eqn:affine-expression-rank-2} we may only see the action of the symmetric group $\Sigma_{2k+1}$ on $(\infty, x_{1}, \ldots, x_{2k+1})$ rather than the full $\Sigma_{2k+2}$. Given a configuration $y = (y_{0}, \ldots, y_{2k+1}) \in \mathrm{Conf}_{2k+2}(\mathbf{P}^{1})$ choose an affine chart $\mathbf{A}^{1} \subset \mathbf{P}^{1}$ so that all $2k+2$ points $y_{0}, \ldots, y_{2k+1}$ lie on it. In a similar way (by replacing $\mathrm{Aff}_{1}(\mathbf{R})$ by $\mathrm{PGL}_{2}(\mathbf{R})$) that we went from \eqref{eqn:gauge-fixed-rank-2} to \eqref{eqn:affine-expression-rank-2} in the proof of \cref{thm:Gelfand-MacPherson-polylogarithmic-form-rank-2} one may see that
  \[
    \widehat{P}_{2}^{2k} = -\frac{4}{(2\pi)^{2k}}\sum_{i < j}\frac{(y_{i} - y_{j})^{2k-2}\log|y_{i}-y_{j}|}{\prod_{\ell\neq i,j}(y_{i}-y_{\ell})(y_{j}-y_{\ell})}\nu_{\mathrm{proj}},
  \]
  where $\nu_{\mathrm{proj}} = \iota_{E_1}\iota_{E_{0}}\iota_{E_{-1}}\Omega$ and
  \[
    E_{-1} = \sum_{i}\partial_{i}, \quad E_{0}=\sum_{i}y_{i}\partial_{i}, \quad E_{1}=\sum_{i}y_{i}^{2}\partial_{i},
  \]
  are the infinitesimal generators of the action of $\mathrm{PGL}_{2}(\mathbf{R})$ in $\mathrm{Conf}_{2k+2}(\mathbf{R})$ and
  \[
    \Omega = \mathrm{d}y_{0} \wedge \cdots \wedge \mathrm{d}y_{2k+1}.
  \]
  In coordinates,
  \[
    \nu_{\mathrm{proj}} = \sum_{a < b < c}(-1)^{a+b+c-3}(y_{b}-y_{a})(y_{c}-y_{a})(y_{c}-y_{b})\mathrm{d}y_{\widehat{abc}},
  \]
  where
  \[
    \mathrm{d}y_{\widehat{abc}} = \mathrm{d}y_{0} \wedge \cdots \wedge \widehat{\mathrm{d}y_{a}} \wedge \cdots \wedge \widehat{\mathrm{d}y_{b}} \wedge \cdots \wedge \widehat{\mathrm{d}y_{c}} \wedge \cdots \wedge \mathrm{d}y_{2k+1}.
  \]
\end{remark}
\begin{remark}\label{rem:Pirio-webs}
  The formula \eqref{eqn:gauge-fixed-rank-2} was conjectured by Pirio, who verified it for $k\leq 6$ \Cite[Proposition 5.28 and the discussion following it]{pirioN+3websRationalCurves2023}. Our theorem proves it for every $k$. To explain its meaning in web geometry, consider the web on $\mathcal{M}_{0,2k+3}(\mathbf{R})$ whose foliations are the fibres of the forgetful maps $\beta_{j}\colon\mathcal{M}_{0,2k+3}(\mathbf{R})\to\mathcal{M}_{0,2k+2}(\mathbf{R})$. By \cref{cor:functional-equation-l-power-first-pontryagin-class}, the signed pullbacks $(-1)^{j}\beta_{j}^{*}\widehat{P}_{2}^{2k}$ sum to zero. They form the Euler abelian relation studied by Damiano \cite{damianoWEBSABELIANEQUATIONS1980,damianoWebsCharacteristicForms1983} and Pirio. As explained in \Cite[Theorem E, part 1]{pirioN+3websRationalCurves2023}, this relation, together with the combinatorial abelian relations coming from subwebs obtained by omitting two foliations, spans the full space of abelian relations. Therefore, our formula makes the remaining relation explicit for every $k$.
\end{remark}
\subsection{Rank $2$ transgression}\label{sec:rank-two-transgression}
Let $P = p_{1}$ the first Pontryagin class. In \Cite{gelfandGeometryGrassmanniansGeneralization1982}, Gelfand and MacPherson obtain a function from $\widehat{P}_{2}^{2}$ by taking primitives, i.e.\ $\frac{2}{\pi^{2}}L_{\mathbf{R}}$ is the primitive of $\widehat{P}_{2}^{2}$. Now, let $P=p_{1}^{k}$ and $k>1$. We need a more involved transgression. Consider the double complex $\Omega^{*}(\mathcal{M}_{0,\bullet}(\mathbf{R}))$ equipped with the de Rham differential $\mathrm{d}\colon \Omega^{*}(\mathcal{M}_{0, \bullet}(\mathbf{R})) \to \Omega^{*+1}(\mathcal{M}_{0, \bullet}(\mathbf{R}))$ and the simplicial differential $\delta\colon \Omega^{*}(\mathcal{M}_{0,\bullet}(\mathbf{R})) \to \Omega^{*}(\mathcal{M}_{0, \bullet + 1}(\mathbf{R}))$ given by the formula $\delta = \sum_{j}(-1)^{j}\beta_{j}^{*}$. Since $\widehat{P}_{2}^{2k}$ is a top form, then $\mathrm{d}\widehat{P}_{2}^{2k} = 0$ and by \cref{cor:functional-equation-l-power-first-pontryagin-class} we have $\delta \widehat{P}_{2}^{2k} = 0$.
The transgression consists of successively taking primitives and applying $\delta$. Since $\mathrm{d}\delta = \delta\mathrm{d}$ and $\delta^{2}=0$, each application of $\delta$ gives a closed form. We seek primitives until we reach a function $\Phi_{0}$, whose simplicial differential is locally constant. Schematically, we want:
\[
  \begin{tikzcd}[column sep=small, row sep=small]
    \widehat{P}_{2}^{2k} \ar[r, "\delta"] & 0 & & & & \\
    \Phi_{2k-2} \ar[u, "\mathrm{d}"] \ar[r, "\delta"] & \delta\Phi_{2k-2} & & & & \\
    & \Phi_{2k-3} \ar[u, "\mathrm{d}"] \ar[r, "\delta"] & \ddots & & & \\
    & & \Phi_{1} \ar[u, "\mathrm{d}"] \ar[r, "\delta"] & \delta\Phi_{1} & & \\
    & & & \Phi_{0} \ar[u, "\mathrm{d}"] \ar[r, "\delta"] & \delta\Phi_{0} \ar[r, "\mathrm{d}"] & 0.
  \end{tikzcd}
\]
Primitives exist in every connected component (also referred as chambers) of $\mathcal{M}_{0, \bullet}(\mathbf{R})$, because each one is contractible. However, it turns out that an application of Brown's work \cite{brownMultipleZetaValues2009} will allow us to show that primitives with polylogarithmic coefficients always exist.
\begin{proposition}\label{prop:Brown-rank-two-transgression}
  For $P=p_{1}^{k}$, there exist forms
  \[
    \Phi_{j} \in \Omega^{j}(\mathcal{M}_{0,4k-j}(\mathbf{R})), \quad 0 \leq j \leq 2k-2,
  \]
  whose coefficients on each chamber are multiple polylogarithms with real constants, such that
  \begin{equation}\label{eqn:rank-two-descent}
    \mathrm{d}\Phi_{2k-2} = \widehat{P}_{2}^{2k}, \quad \mathrm{d}\Phi_{j} = \delta\Phi_{j+1} \quad (0 \leq j < 2k-2).
  \end{equation}
  In particular, $\delta\Phi_{0}$ is locally constant on $\mathcal{M}_{0,4k+1}(\mathbf{R})$.
\end{proposition}
\begin{proof}
  Brown's algebra of multiple polylogarithms on $\mathcal{M}_{0,n}$, with a choice of branch, has vanishing de Rham cohomology in positive degrees \Cite[Theorem 6.26]{brownMultipleZetaValues2009}. Thus closed forms with coefficients in this algebra admit polylogarithmic primitives. We use its real branches on each chamber and extend the constants to $\mathbf{R}$. This class of forms is preserved by pullbacks under the forgetful maps: the logarithmic one-forms defining the iterated integrals pull back to logarithmic one-forms on the larger configuration space.

  Formula \eqref{eqn:affine-expression-rank-2} belongs to this complex. Indeed, on a real chamber every $x_{i}-x_{j}$ has constant sign, so that $\log|x_{i}-x_{j}| = \log(\varepsilon_{ij}(x_{i}-x_{j}))$ for a fixed $\varepsilon_{ij}\in\{1,-1\}$. Since $\widehat{P}_{2}^{2k}$ is closed, Brown's theorem gives $\Phi_{2k-2}$. Then
  \[
    \mathrm{d}\delta\Phi_{2k-2} = \delta\mathrm{d}\Phi_{2k-2} = \delta\widehat{P}_{2}^{2k} = 0,
  \]
  and Brown's theorem applies again. At each subsequent step, $\mathrm{d}\delta\Phi_{j} = \delta^{2}\Phi_{j+1}=0$, so we obtain \eqref{eqn:rank-two-descent}. Finally, $\mathrm{d}\delta\Phi_{0}=0$, which proves the last assertion.
\end{proof}
This argument does not imply that only dilogarithms occur. Nevertheless, for $k=2$ the entire descent can be expressed using Rogers dilogarithms and logarithms, as follows. We work on the chamber
\[
  (\infty,0,1,x_{3},\ldots,x_{n-1}), \quad 1 < x_{3} < \cdots < x_{n-1}.
\]
On $\mathcal{M}_{0,6}(\mathbf{R})$, set
\begin{equation}\label{eqn:rank-two-k2-first-primitive}
  \Phi_{2} = \frac{1}{240\pi^{4}}\operatorname{Alt}_{6}\left[L_{\mathbf{R}}(x_{3})\mathrm{d}\log|x_{4}|\wedge\mathrm{d}\log|x_{5}|\right],
\end{equation}
where $\operatorname{Alt}_{6}$ is the unnormalised signed sum over permutations of the six marked points, renormalising the first three to $(\infty,0,1)$ after each permutation. On $\mathcal{M}_{0,7}(\mathbf{R})$, set
\[
  \Phi_{1} = \frac{1}{120\pi^{2}}\sum_{(f,g)}\left(\log f\,\mathrm{d}\log g - \log g\,\mathrm{d}\log f\right),
\]
where the sum runs over the seven pairs
\[
  \begin{gathered}
    \left(\frac{x_{3}x_{5}(x_{4}-1)}{x_{4}(x_{5}-x_{3})},x_{6}\right), \quad
    \left(\frac{x_{4}(x_{5}-x_{3})}{x_{5}(x_{4}-x_{3})},x_{6}-x_{3}\right), \\
    \left(\frac{x_{4}(x_{4}-x_{3})}{x_{3}(x_{4}-1)},x_{6}-x_{4}\right), \quad
    \left(\frac{1}{x_{4}},x_{6}-1\right), \\
    \left(\frac{x_{4}-1}{x_{3}-1},x_{5}-1\right), \quad
    \left(\frac{x_{3}-1}{x_{3}},x_{5}-x_{3}\right), \quad (x_{3},x_{4}).
  \end{gathered}
\]
All these functions are positive on the chosen chamber. The terminal function on $\mathcal{M}_{0,8}(\mathbf{R})$ is
\[
  \begin{aligned}
    \Phi_{0} = \frac{1}{60\pi^{2}}\bigg[&-2L_{\mathbf{R}}(x_{3}) + L_{\mathbf{R}}(x_{4}) - L_{\mathbf{R}}(x_{5}) \\
    &-L_{\mathbf{R}}(x_{4}/x_{3}) + L_{\mathbf{R}}(x_{5}/x_{3}) - L_{\mathbf{R}}(x_{5}/x_{4})\bigg].
  \end{aligned}
\]
These expressions satisfy
\begin{equation}\label{eqn:rank-two-k2-descent}
  \mathrm{d}\Phi_{2} = \widehat{P}_{2}^{4}, \quad \mathrm{d}\Phi_{1} = \delta\Phi_{2}, \quad \mathrm{d}\Phi_{0} = \delta\Phi_{1}, \quad \delta\Phi_{0} = \frac{1}{120}.
\end{equation}
Indeed, differentiating \eqref{eqn:rank-two-k2-first-primitive} gives \eqref{eqn:affine-expression-rank-2} for $k=2$. Expanding the forgetful pullbacks and using the Rogers five-term relation \eqref{eqn:5-term-relation} gives
\[
  \delta\Phi_{2} = \frac{1}{60\pi^{2}}\sum_{(f,g)}\mathrm{d}\log f\wedge\mathrm{d}\log g = \mathrm{d}\Phi_{1}.
\]
The next identity follows by expanding $\delta\Phi_{1}$ and differentiating \eqref{eqn:rogers-dilog}, which gives
\[
  \mathrm{d}L_{\mathbf{R}}(u) = \frac{1}{2}\left(\log|u|\mathrm{d}\log|1-u| - \log|1-u|\mathrm{d}\log|u|\right).
\]
To check the last identity, write the five-term expression in \eqref{eqn:5-term-relation}, for $\Phi=L_{\mathbf{R}}$, as
\[
  \begin{aligned}
    E(u,v) = {}&-L_{\mathbf{R}}\left(\frac{u(v-1)}{v(u-1)}\right) + L_{\mathbf{R}}\left(\frac{v-1}{u-1}\right) - L_{\mathbf{R}}(v/u) \\
    &+L_{\mathbf{R}}(v)-L_{\mathbf{R}}(u).
  \end{aligned}
\]
For $1<u<v$ we have $E(u,v)=-\pi^{2}/6$, and direct cancellation gives
\[
  \begin{aligned}
    60\pi^{2}\delta\Phi_{0} = {}&-2E(x_{3},x_{4})+E(x_{3},x_{5})-E(x_{3},x_{6}) \\
    &-E(x_{4},x_{5})+E(x_{4},x_{6})-E(x_{5},x_{6}).
  \end{aligned}
\]
Since the coefficients sum to $-3$, we obtain $\delta\Phi_{0}=1/120$. The endpoint can be simplified further by changing the intermediate primitive. The expression
\[
  H = \Phi_{0}+\frac{1}{20\pi^{2}}L_{\mathbf{R}}(x_{3})
\]
depends only on the first six marked points, so we may regard it as a function on $\mathcal{M}_{0,7}(\mathbf{R})$. Its coboundary when regarded as a function of eight points vanishes, since $\delta\Phi_{0}=1/120$ and $\delta L_{\mathbf{R}}(x_{3})=-\pi^{2}/6$. Cancelling the final deletion terms therefore gives, on $\mathcal{M}_{0,8}(\mathbf{R})$,
\[
  \delta H = \Phi_{0}+\frac{1}{20\pi^{2}}L_{\mathbf{R}}(x_{3}).
\]
Consequently, the replacements
\[
  \Phi_{1}' = \Phi_{1}-\mathrm{d}H, \quad \Phi_{0}' = \Phi_{0}-\delta H = -\frac{1}{20\pi^{2}}L_{\mathbf{R}}(x_{3})
\]
preserve \eqref{eqn:rank-two-k2-descent}. If $\Phi = \frac{2}{\pi^{2}}L_{\mathbf{R}}$ denotes the $k=1$ primitive from \cref{rem:specialisation-rogers-dilog}, then $\delta\Phi=-1/3$ on these chambers and
\[
  \Phi_{0}' = \frac{3}{40}\Phi\smile\delta\Phi.
\]
Here we use the usual cup product of cochains: for cochains of degrees $p$ and $q$, respectively,
\[
  (A\smile B)(z_{0},\ldots,z_{p+q}) = A(z_{0},\ldots,z_{p})B(z_{p},\ldots,z_{p+q}).
\]
In particular, $\Phi$ has cochain degree $3$ and $\delta\Phi$ has degree $4$.
\begin{question}\label{conj:rank-two-Rogers-descent}
  For every $k\geq 1$, can the descent \eqref{eqn:rank-two-descent} be chosen with terms given by finite expressions in Rogers dilogarithms of rational configuration functions, logarithms, logarithmic differentials? Moreover, on the compatible ordered chambers, can the primitives be chosen so that
  \[
    \Phi_{0} = c_{k}\Phi\smile(\delta\Phi)^{\smile(k-1)} + C,
  \]
  where $\Phi = \frac{2}{\pi^{2}}L_{\mathbf{R}}$, $c_{k}$ depends only on $k$, and $C$ is locally constant? The normalisations above give $c_{1}=1$ and $c_{2}=3/40$.
\end{question}
The cup-product identity
\[
  \delta\left(\Phi\smile(\delta\Phi)^{\smile(k-1)}\right) = (\delta\Phi)^{\smile k}
\]
already gives a $(4k+1)$-term relation derived from the Rogers five-term relation \eqref{eqn:5-term-relation}. The question asks whether this expression, up to a scalar and a locally constant term, can occur as the endpoint of a descent of the specific Gelfand-MacPherson form $\widehat{P}_{2}^{2k}$. We leave questions of canonicity and uniqueness open.
\subsection{Gelfand-MacPherson forms for self-dual Grassmannians and the Pontryagin character}
We normalise the Pontryagin character by
\[
  \mathrm{ph}_{k}(E) = \mathrm{ch}_{2k}(E\otimes_{\mathbf{R}}\mathbf{C}).
\]
Newton's identities give
\begin{equation}\label{eqn:newton-identity}
  \mathrm{ph}_{k} = \frac{(-1)^{k-1}}{(2k-1)!}p_{k} + \text{decomposable Pontryagin classes}.
\end{equation}
Although the corresponding characteristic forms differ by these decomposable terms, their Gelfand-MacPherson one-forms are proportional.
\begin{proposition}\label{prop:GM-character-Pontryagin-comparison}
  If $p+q=4k$, then
  \[
    \widehat{(\mathrm{ph}_{k})}_{p}^{q} = \frac{(-1)^{k-1}}{(2k-1)!}\widehat{(p_{k})}_{p}^{q}.
  \]
\end{proposition}
\begin{proof}
  Put $m=p+q-1=4k-1$. At a point $\xi$ of the Grassmannian, let $V_{\xi}$ be the tangent space to the Cartan orbit. Consider two Pontryagin forms $\alpha$ and $\beta$ (each of positive degree), with $\deg\alpha+\deg\beta = 4k$. By \cref{thm:vanishing-orbits}, $\alpha|_{V_{\xi}} = 0$ and $\beta|_{V_{\xi}} = 0$. We claim that
  \[
    (\alpha\wedge\beta)_{\xi}(v_{1},\ldots,v_{m},w) = 0
  \]
  for every $v_{1},\ldots,v_{m}\in V_{\xi}$ and every ambient tangent vector $w$. Indeed, contraction satisfies
  \[
    \iota_{w}(\alpha\wedge\beta) = (\iota_{w}\alpha)\wedge\beta + (-1)^{\deg\alpha}\alpha\wedge(\iota_{w}\beta).
  \]
  Restricting this identity to $V_{\xi}$ proves the claim. Therefore, fibre integration along $\pi$ for $\alpha \wedge \beta$ vanishes. Every decomposable Pontryagin monomial can be written as such a product. Thus, fibre integration kills all decomposable terms in Newton's identity \eqref{eqn:newton-identity}. Applying $\pi_{*}$ gives the desired identity.
\end{proof}
Damiano's vanishing theorem gives a stronger restriction on the Pontryagin character.
\begin{theorem}[{\Cite{damianoVanishingCharacteristicForms1982}}]\label{thm:damiano-character-vanishing}
  Let $P=\mathrm{ph}_{k}$, where $k\geq 1$, and let $p,q\geq 1$ with $p+q-1\leq 4k$. If
  \[
    4k-(p+q-1) < \frac{4k}{1+2\min(p,q)},
  \]
  then $P_{p}^{q}$ vanishes when evaluated on $p+q-1$ vectors tangent to a Cartan orbit and arbitrary remaining tangent vectors. In particular, $\widehat{P}_{p}^{q}=0$.
\end{theorem}
\begin{remark}\label{rem:numerology-one-form}
  For $p+q=4k$, this inequality holds whenever $p\neq q$. Thus the self-dual Grassmannian $\mathrm{G}_{2k}^{2k}(\mathbf{R})$ is the only Grassmannian on which the associated Gelfand-MacPherson one-form can be nonzero.
\end{remark}
\begin{corollary}\label{cor:self-dual-GM-functional-equations}
  For $P=\mathrm{ph}_{k}$ or $P=p_{k}$, the form $\widehat{P}_{2k}^{2k}$ satisfies the dual functional equations
  \[
    \sum_{j=1}^{4k+1}(-1)^{j}\alpha_{j}^{*}\widehat{P}_{2k}^{2k} = 0 \quad \text{on } \mathrm{C}_{2k+1}^{2k}(\mathbf{R}),
  \]
  and
  \[
    \sum_{j=1}^{4k+1}(-1)^{j}\beta_{j}^{*}\widehat{P}_{2k}^{2k} = 0 \quad \text{on } \mathrm{C}_{2k}^{2k+1}(\mathbf{R}).
  \]
\end{corollary}
\begin{proof}
  By \cref{cor:0-form-vanishing}, the Gelfand-MacPherson $0$-forms on $\mathrm{C}_{2k+1}^{2k}(\mathbf{R})$ and $\mathrm{C}_{2k}^{2k+1}(\mathbf{R})$ vanish. Apply the differential relation \eqref{eqn:differential-relation} to these forms for $P=\mathrm{ph}_{k}$. Applying \cref{thm:damiano-character-vanishing} gives the pair of functional equations. The comparison in \cref{prop:GM-character-Pontryagin-comparison} transfers both identities to $p_{k}$.
\end{proof}
We now compute the integrand. Represent a configuration by
\[
  S(x) = (\mathrm{I}_{2k\times 2k}\mid\mathbf{1}\mid z_{1}\mid\cdots\mid z_{2k-1}) = (v_{1},\ldots,v_{4k}), \quad z_{b} = (x_{1b},\ldots,x_{2k-1,b},1)^{t}.
\]
Parametrise its positive Cartan orbit by the interior of the simplex
\[
  \Delta^{4k-1} = \left\{r_{i}\geq 0:\sum_{i=1}^{4k}r_{i}=1\right\},
\]
using the row space of
\[
  A(r,x) = S(x)\operatorname{diag}(\sqrt{r_{1}},\ldots,\sqrt{r_{4k}}).
\]
Set
\begin{equation}\label{eqn:self-dual-Gram-determinant}
  G = AA^{t} = \sum_{i=1}^{4k}r_{i}v_{i}v_{i}^{t}, \quad Q(r,x) = \det G = \sum_{|I|=2k}r_{I}\Delta_{I}(S(x))^{2},
\end{equation}
where $r_{I}=\prod_{i\in I}r_{i}$ and the last equality is Cauchy-Binet. We orient the simplex by $\mathrm{d}r_{1}\wedge\cdots\wedge\mathrm{d}r_{4k-1}$, with $r_{4k}=1-\sum_{i<4k}r_{i}$; this agrees with the Cartan fibre orientation. Let $\mathrm{d}$ differentiate in both $r$ and $x$, and introduce the skew-symmetric matrix of two-forms
\[
  \Theta = \mathrm{d}A\wedge(\mathrm{I}_{4k\times 4k}-A^{t}G^{-1}A)\mathrm{d}A^{t}.
\]
Define $N_{ab}(r,x)$ by
\begin{equation}\label{eqn:self-dual-simplex-numerator}
  [Q^{k}\det\Theta]_{(4k-1,1)} = \mathrm{d}r_{1}\wedge\cdots\wedge\mathrm{d}r_{4k-1}\wedge\sum_{a,b=1}^{2k-1}N_{ab}(r,x)\mathrm{d}x_{ab},
\end{equation}
where the subscript selects the component of fibre degree $4k-1$ and base degree $1$. Products of differential forms in the determinant are wedge products. The following theorem shows, in particular, that the coefficients $N_{ab}$ are polynomial.
\begin{theorem}\label{thm:self-dual-GM-simplex-integrand}
  With the numerator defined by \eqref{eqn:self-dual-simplex-numerator}, the Gelfand-MacPherson form for $P = \mathrm{ph}_{k}$ is
  \[
    \widehat{P}_{2k}^{2k} = \frac{(-1)^{k-1}2^{4k-1}}{(2k-1)!(2\pi)^{2k}}\sum_{a,b=1}^{2k-1}\left(\int_{\Delta^{4k-1}}\frac{N_{ab}(r,x)}{Q(r,x)^{k+1}}\mathrm{d}r_{1}\cdots\mathrm{d}r_{4k-1}\right)\mathrm{d}x_{ab},
  \]
  where $N_{ab}(r,x)$ are polynomials in $r$ and $x$.
\end{theorem}
\begin{proof}
  In the orthonormal frame obtained from $A$, the tautological connection has curvature
  \[
    \Omega = G^{-1/2}\Theta G^{-1/2}.
  \]
  On the oriented double cover, the top Pontryagin form of the rank-$2k$ tautological bundle is the square of its Euler form. Hence
  \[
    (p_{k})_{2k}^{2k} = \frac{\det\Omega}{(2\pi)^{2k}} = \frac{\det\Theta}{(2\pi)^{2k}Q}.
  \]
  We show that $Q^{k}\det\Theta$ has polynomial coefficients before imposing $\sum r_{i}=1$. The square roots cancel in $\mathrm{d}AA^{t}$ and $\mathrm{d}A\wedge\mathrm{d}A^{t}$. Since
  \[
    \Theta = \mathrm{d}A\wedge\mathrm{d}A^{t} - (\mathrm{d}AA^{t})G^{-1}\wedge(A\mathrm{d}A^{t}), \qquad G^{-1} = \frac{\operatorname{adj}(G)}{Q},
  \]
  the coefficients of $\Theta$ have denominators only among powers of $Q$. Fix $i$ and write
  \[
    G = B+r_{i}v_{i}v_{i}^{t}, \qquad B = \sum_{\ell\ne i}r_{\ell}v_{\ell}v_{\ell}^{t}.
  \]
  We temporarily allow rational functions of all variables except $r_{i}$. The remaining columns span $\mathbf{R}^{2k}$, so $\det B$ is not the zero polynomial. Thus $B^{-1}$ is allowed and introduces no denominators involving $r_{i}$. Direct matrix multiplication gives
  \[
    G^{-1} = B^{-1} - \frac{r_{i}\det B}{Q}(B^{-1}v_{i})(B^{-1}v_{i})^{t}.
  \]
  Hence,
  \[
    \Theta_{ab} = \text{a two-form polynomial in }r_{i} + \frac{r_{i}\det B}{Q}\eta_{a}\wedge\eta_{b}, \qquad \eta = (\mathrm{d}AA^{t})B^{-1}v_{i}.
  \]
  The coefficients of each $\eta_{a}$ are also polynomial in $r_{i}$. Now expand the $2k\times 2k$ determinant. A term using $j$ fractions contains $2j$ one-forms chosen from $\eta_{1},\ldots,\eta_{2k}$. If $j>k$, one of these forms occurs twice, so the term vanishes. Every surviving term therefore has denominator at most $Q^{k}$. It follows that $Q^{k}\det\Theta$ is polynomial in $r_{i}$.

  The same argument applies to every $r_{i}$. Thus any denominator remaining in a reduced coefficient of $Q^{k}\det\Theta$ could depend only on $x$. But this denominator must also divide a power of $Q$. No nonconstant polynomial depending only on $x$ divides $Q$, because the coefficient of $r_{1}\cdots r_{2k}$ in $Q$ is $1$. Therefore $Q^{k}\det\Theta$ has polynomial coefficients in $r$ and $x$. Substituting $r_{4k}=1-\sum_{i<4k}r_{i}$ and taking the component in \eqref{eqn:self-dual-simplex-numerator} proves that each $N_{ab}$ is polynomial. It follows that the component contributing to fibre integration is
  \[
    \frac{\mathrm{d}r_{1}\wedge\cdots\wedge\mathrm{d}r_{4k-1}\wedge\sum_{a,b}N_{ab}\mathrm{d}x_{ab}}{(2\pi)^{2k}Q^{k+1}}.
  \]
  Integrating this component over the simplex, multiplying by $2^{4k-1}$ and applying \cref{prop:GM-character-Pontryagin-comparison} proves the formula.
\end{proof}
\begin{remark}\label{rem:self-dual-Gram-quadric}
  The denominator in \cref{thm:self-dual-GM-simplex-integrand} has the same origin as the Gram determinant \eqref{eqn:Gram-determinant-rank-2} in the rank-two calculation. In the Pl\"ucker coordinates $\Delta_{I}$ of the rescaled matrix $A$, Cauchy-Binet gives
  \[
    Q = \sum_{|I|=2k}\Delta_{I}^{2}.
  \]
  Thus, the possible poles lie over the \textsf{Gram quadric}, whose intersection with the complex Grassmannian parametrises planes on which the standard symmetric bilinear form becomes degenerate. Its pullback to the simplex is the degree-$2k$ polynomial \eqref{eqn:self-dual-Gram-determinant}. For rank two this is quadratic in $r$. In higher rank it remains quadratic in Pl\"ucker coordinates.
\end{remark}

\subsection{Closedness of the Gelfand-MacPherson forms in self-dual Grassmannians}\label{subsec:closedness}
In \Cite[\textsection 1.2]{goncharovSimpleConstructionGrassmannian2013}, Goncharov describes Gelfand and MacPherson's construction as a real-valued Grassmannian $2k$-logarithm obtained by constructing its differential. Every homogeneous Pontryagin polynomial of degree $4k$ is a multiple of $\mathrm{ph}_{k}$ modulo decomposable terms. By the orbit-vanishing argument in the proof of \cref{prop:GM-character-Pontryagin-comparison}, these decomposable terms have zero fibre integral whenever $p+q=4k$. Thus every Gelfand-MacPherson one-form arising from such a polynomial is a scalar multiple of the one associated to $\mathrm{ph}_{k}$. Moreover, \cref{thm:damiano-character-vanishing} forces this one-form to vanish unless $p=q=2k$. Therefore, the only possible nonzero one-form to which Goncharov's description can refer is, up to a scalar, $\widehat{P}_{2k}^{2k}$ for $P=\mathrm{ph}_{k}$. This description does not even define a local function unless the one-form is shown to be closed. The known vanishing theorems do not establish this closedness in general. The obstruction is given by the following proposition.
\begin{proposition}\label{prop:self-dual-closedness-obstruction}
  For $P = \mathrm{ph}_{k}$, the $k$th component of the Pontryagin character, we have that
  \[
    \mathrm{d}\widehat{P}_{2k}^{2k} = \frac{2(-1)^{k-1}}{(2k-1)!}\sum_{j=1}^{4k}(-1)^{j}\beta_{j}^{*}\left(\sum_{a=1}^{k-1}\widehat{(p_{a}p_{k-a})}_{2k}^{2k-1}\right).
  \]
\end{proposition}
\begin{proof}
  The tautological bundle on $\mathrm{G}_{2k-1}^{2k}(\mathbf{R})$ has rank $2k-1$, so its $p_{k}$-form vanishes. Applying \eqref{eqn:differential-relation} to $p_{k}$ gives
  \[
    \mathrm{d}\widehat{(p_{k})}_{2k}^{2k} = 2\sum_{j=1}^{4k}(-1)^{j}\beta_{j}^{*}\widehat{(p_{k})}_{2k}^{2k-1}.
  \]
  On $\mathrm{G}_{2k}^{2k-1}(\mathbf{R})$, let $E:=\gamma_{2k}^{2k-1}$ be the tautological bundle. Its orthogonal complement has rank $2k-1$. Hence
  \[
    p(E)p(E^{\perp}) = 1, \qquad p_{k}(E^{\perp}) = 0.
  \]
  Taking the degree-$4k$ component gives
  \[
    p_{k}(E) = \sum_{a=1}^{k-1}p_{a}(E)p_{k-a}(E) + \text{products of at least three positive-degree Pontryagin forms}.
  \]
  Fibre integration on $\mathrm{G}_{2k}^{2k-1}(\mathbf{R})$ evaluates these forms on $4k-2$ vertical vectors and two remaining ambient tangent vectors. In every term with at least three factors, one factor will always have only vertical vectors as inputs. Hence, that factor vanishes by \cref{thm:vanishing-orbits}. Therefore, only the quadratic factors contribute. The comparison in \cref{prop:GM-character-Pontryagin-comparison} now gives the assertion.
\end{proof}
\begin{remark}\label{rem:k-one-closedness}
For $k=1$, the sum is empty and \cref{prop:self-dual-closedness-obstruction} recovers closedness. For $k=2$, it gives
\[
  \mathrm{d}\widehat{(\mathrm{ph}_{2})}_{4}^{4} = -\frac{1}{3}\sum_{j=1}^{8}(-1)^{j}\beta_{j}^{*}\widehat{(p_{1}^{2})}_{4}^{3}.
\]
The orbit-vanishing argument does not eliminate these terms: each of the two $p_{1}$-factors can receive one nonvertical vector. Damiano's vanishing in \cref{thm:damiano-character-vanishing} does not establish the vanishing of this signed sum.
\end{remark}

Thus closedness in higher rank requires a vanishing identity beyond those supplied by \cref{thm:vanishing-orbits,thm:damiano-character-vanishing}. Without establishing this identity, the interpretation of these one-forms as differentials of ``Gelfand-MacPherson Grassmannian polylogarithms'' remains unjustified. We leave this as a question.
\begin{question}\label{question:self-dual-GM-closedness}
  For every $k\geq 2$, is the Gelfand-MacPherson one-form $\widehat{P}_{2k}^{2k}$ for $P=\mathrm{ph}_{k}$ closed? Equivalently, does the alternating sum in \cref{prop:self-dual-closedness-obstruction} vanish?
\end{question}
\printbibliography
\end{document}